\documentclass[fleqn,10pt]{wlscirep}
\pdfoutput=1
\usepackage{subfigure}  
\usepackage{color}
\usepackage{mathtools}
\usepackage{amssymb}
\usepackage{graphicx}

\title{Correcting Observation Model Error\\ in Data Assimilation}
\author[1]{Franz Hamilton}
\author[2]{Tyrus Berry}
\author[2,*]{Timothy Sauer}
\affil[1]{North Carolina State University, Department of Mathematics, Raleigh, 27695, USA}
\affil[2]{George Mason University, Department of Mathematical Sciences, Fairfax, 22030, USA}

\affil[*]{tsauer@gmu.edu}

\begin{abstract}
Standard methods of data assimilation assume prior knowledge of a model that describes the system dynamics and an observation function that maps the model state to a predicted output. An accurate mapping from model state to observation space is crucial in filtering schemes when adjusting the estimate of the system state during the filter's analysis step. However, in many applications the true observation function may be unknown and the available observation model may have significant errors, resulting in a suboptimal state estimate. We propose a method for observation model error correction within the filtering framework. The procedure involves an alternating minimization algorithm used to iteratively update a given observation function to increase consistency with the model and prior observations, using ideas from attractor reconstruction. The method is demonstrated on the Lorenz 1963 and Lorenz 1996 models, and on a single-column radiative transfer model with multicloud parameterization.
\end{abstract}
\begin{document}

\flushbottom
\maketitle

\thispagestyle{empty}

\section{Introduction}
Data assimilation plays an increasingly important role in nonlinear science, as a means of inferring unobserved model variables and constraining unknown parameters. Use of the Extended Kalman Filter (EKF) and Ensemble Kalman Filter (EnKF) is now standard in a wide range of geophysical problems \cite{enkf7,evensen,julier1,julier2,rabier,TELA:TELA066,cummings} and several areas of physical and biological sciences where spatiotemporal dynamics is involved \cite{yoshida,stuart,schiffbook,berry2}.


 These standard applications require complete knowledge of the system equations and observation functions. Current research investigates the effects of incomplete knowledge on this process, such as model error, missing equations and multiple sources of error in observations. In particular, the issue of observation errors, due to truncation, resolution differences, and instrument error, has received great recent attention \cite{dee,Hodyss1,Hodyss2,van2015representation,JanjicRepresentation,BerrySauer18}. In the case of unknown or incorrect observation models, there is interest in fixing these deficiencies.
 For example, a recent study \cite{berry2017} discusses replacing an unknown observation function with a training set of observations and accompanying states.

In this article, an iterative approach to fixing observation model error is proposed which does not require training data, and can be applied as part of a sequential data assimilation implementation. The idea is based on an alternating minimization algorithm applied to the observation function.  In the first step, a filter (eg. Kalman-type or variational filters) is applied to find the optimal state estimate based on the given observation model.  In the second step, an observation model correction term is interpolated from the difference between the actual observations and the observation model applied to the state estimate produced by the filter; this interpolation is localized in the underlying phase space of the dynamical system.  The model correction term is then applied to form a new observation model. The two steps are then repeated until convergence.

Fig.~\ref{figure3} shows an example application of the technique, to the Lorenz attractor with dynamical noise. The underlying  model equations (the Lorenz equations) are assumed known. An initial guess is made for the observation function used in the filter, which is far from the function generating the observed data. Sequential filtering is applied iteratively, and the observation model correction is learned through the iteration. The RMSE of the filter decreased with iteration number, and after about a dozen iterations the minimum RMSE is approximately attained.

Several other examples illustrate the varying contexts in which the method can be applied. A critical hurdle for all filtering methods is the ability to scale up to large problems, which is typically achieved with a spatial localization.  As a test case for spatiotemporal data we consider the Lorenz-96 system, in networks with 10 and 40 nodes. In the latter case, a spatial localization technique is developed which allows interpolation within each local region.  Finally, we consider a more physically realistic example where observation model error can be especially detrimental to filtering, namely the case of radiative transfer models (RTM). To simulate severe observation model error, we assign the cloud fractions of a typical RTM to zero in the observation model.  We then generate data using the full RTM (including the cloud fractions) and apply our method using the crippled observation function (with cloud fractions set to zero). The results show significant improvement in RMSE after three iterations of our observation model error correction algorithm.

The algorithm for correcting the observation model error is described in Section 2, along with its relation to alternating minimization methods in optimization theory, and details of its implementation in an ensemble filter. Sections 3 and 4 describe applications of the algorithm to Lorenz-63 and Lorenz-96 models, the latter to show how the method scales for spatiotemporal problems. The application to the radiative transfer model in shown in Section 5.

\section{Filtering with an incorrect observation function}

In the general filtering problem, we assume a system with $n$-dimensional state vector $x$ and $m$-dimensional observation vector $y$ defined by 
\begin{eqnarray} \label{e1}
x_{k} &=& f(x_{k-1})+w_{k-1}\nonumber \\
y_{k} &=& h(x_{k})+v_{k}
\end{eqnarray}
where $w_{k-1}$ and $v_{k}$ are white noise processes with covariance matrices $Q$ and $R$, respectively. The function $f$ represents the system dynamics and $h$ is an observation function that maps the model state to a predicted output. The goal is to sequentially estimate the state of the system given some noisy observations.  Below we will consider a specific filtering algorithm, however, at this point our approach can be formulated in terms of a generic filtering method.

\subsection{The observation error correction algorithm}



The effectiveness of standard filtering approaches is based on the assumption that the observation function $h$ is perfectly known. The goal of this section is to address what happens when $h$ is {\it not known}, and in its place an incorrect observation function $g$ is used. In fact, observation model errors can have many sources, from truncation error due to downsampling high resolution state variables (also called representation error) to simple mismatch between the actual and available observation functions (often referred to as observation model error) \cite{Hodyss1,van2015representation,JanjicRepresentation}.  In this article we will take a very general outlook by considering $h$ to be the true mapping from the fully resolved true state variables $x_k$ into observed variables $y_k$, which is subject only to instrument error $v_k$.  Meanwhile, $g$ will denote a possibly incorrect mapping from state variables into observation variables which can be compared to the actual observations $y_k$.  In such a situation, we can rewrite of the second part of Eq. \ref{e1} as
\begin{eqnarray}
\label{e2}
y_{k} &=& h(x_{k}) + v_{k} \nonumber\\
&=& g(x_{k})+b(x_k) + v_{k}
\end{eqnarray}
where $b$ is the error in our estimate resulting from use of the incorrect observation function.  The term $b(\cdot)$ encapsulates all sources of error except for instrument noise which is the noise term $v_k$. We can write this error term as $b(x_k) = h(x_k)-g(x_k)$, or the difference between the true and incorrect observation functions at step $k$. Note that this error is dependent on the fully resolved state $x_k$.

Repairing observation model error was addressed recently \cite{berry2017} by building a nonparametric estimate of the function $b$ using a training set consisting of observations along with the corresponding true state.  In the current article, we assume that the true state is {\it not} available. A novel approach will be proposed for empirically estimating the model error term $b$ using only the observations $y_k$.  
We begin by describing our method generically for any filtering scheme. The general idea is to iteratively update the incorrect observation function $g$ by obtaining successively improved estimates of the observation model error.


We make an initial definition $g^{(0)}  = g$. The filter is given the known system dynamics $f$, the initial incorrect observation function $g^{(0)}$, and the observations $y$, and provides an estimate of the state at each observation time $k$, which we denote $x_{k}^{(0)}$. This initial state estimate will be subject to large errors, due to the unaccounted-for observation model error. Using this imperfect state estimate, we calculate a noisy estimate  $\hat{b}_{k}^{(0)}$ of the observation model error, corresponding to observation $y_k$ where
\begin{eqnarray}
\label{e3}
\hat{b}_{k}^{(0)} = y_k - g\left(x_{k}^{(0)}\right).
\end{eqnarray}

Due to noise in the data as well as the imperfection of the state estimate, $\hat{b}_{k}^{(0)}$ will not accurately reflect the true observation model error, $b(x_k)$. To build a better estimate of $b(x_k)$, we use a  standard method of nonparametric attractor reconstruction \cite{takens,PCFS,SYC,sauer04} to interpolate the observation model error function, as follows. Given observation $y_k$, we introduce the delay-coordinate vector $z_k = \left[y_k,y_{k-1},\hdots,y_{k-d}\right]$ where $d$ is the number of delays. The delay vector $z_k$ represents the state of the system \cite{takens,SYC}.  To build the reconstruction, we locate the $N$ nearest neighbors (with respect to Euclidean distance) $z_{k_1},...,z_{k_N}$, where
\[ z_{k_j} = [y_{k_j}, y_{k_j-1}, \ldots, y_{k_j-d}] \]
within the set of observations. Once the neighbors are found, the corresponding $\hat{b}_{k_1}^{(0)},\hat{b}_{k_2}^{(0)}, \ldots, \hat{b}_{k_N}^{(0)}$ values are used to estimate $b(x_k)$ by the weighted average
\begin{eqnarray}
\label{e4}
b^{(0)}(x_k) = w_{k_1}\hat{b}_{k_1}^{(0)}+w_{k_2}\hat{b}_{k_2}^{(0)}+\hdots+w_{k_N}\hat{b}_{k_N}^{(0)}.
\end{eqnarray}
which is a locally constant model.  The weights can be chosen according to many strategies, and in order to impose some smoothness on the function $b^{(0)}$, we will use weights which decay exponentially in the distance in delay space, namely, the weight for $j^{th}$ neighbor is defined as
\begin{eqnarray*}
w_{k_j} = \frac{e^{-||z_{k_j}-z_k||/\sigma}}{\sum_{j=1}^{N}e^{-||z_{k_j}-z_k||/\sigma}}.
\end{eqnarray*}
Here, $||z_{k_j}-z_k||$ is the distance of the $j$-th nearest neighbor, $z_{k_j}$, to the current delay-coordinate vector, $z_k$, and $\sigma$ is the bandwidth which controls the weighting of the neighbors in the local model. While many methods are available to tune the $\sigma$ variable, we simply set it to half of the mean distance of the $N$ nearest neighbors to give a smooth roll off of the weights with distance which adapts to naturally to the locally density of the data.

Note that Eq. (\ref{e4}) is still just an approximation of $b(x_k)$, although a more accurate estimate compared to Eq. (\ref{e3}). Our observation function can now be updated, namely
\begin{eqnarray*}
g^{(1)} = g+b^{(0)}.
\end{eqnarray*}
This improved observation function is given to the filter, and the data are re-processed. An improved state estimate, $x^{(1)}_k$, at time $k$ is obtained, a more accurate reconstruction, $b^{(1)}(x_k)$, of the observation model error is formed using Eqs. (\ref{e3}-\ref{e4}) and the observation function is again updated, $g^{(2)} = g+b^{(1)}$.

The method continues iteratively, each iteration an improved reconstruction of $b(x_k)$ is obtained resulting in a better estimate of the state on the next iteration. The method is summarized for steps $\ell = 0, 1, 2, \ldots$ as follows:
\begin{enumerate}
\item Initialize $g^{(0)} = g$, $\Delta g=\textup{Inf}$
\item While $\Delta g$ is greater than threshold
\begin{enumerate}
\item For each observation $y_k$, use filter to estimate state $x^{(\ell)}_k$ given known $f$ and observation function $g^{(\ell)}$
\item Calculate the noisy observation model error estimates $\hat{b}^{(\ell)}_k = y_k-g(x^{(\ell)}_k)$
\item For each $k$, find the $N$-nearest neighbors of delay vector $z_k$ and set
\begin{eqnarray}
\label{e34}
b^{(\ell)}(x_k) &=& w_{k_1}\hat{b}_{k_1}^{(\ell)}+w_{k_2}\hat{b}_{k_2}^{(\ell)}+\hdots+w_{k_N}\hat{b}_{k_N}^{(\ell)}
\end{eqnarray} 
\item Update the observation function, $g^{(\ell+1)} = g+b^{(\ell)}$
\item Update $\Delta g = \frac{1}{T}\sum_{k=1}^T |\hat b^{(\ell)}_k - \hat b^{(\ell-1)}_k|$
\end{enumerate}
\end{enumerate}

In the absence of results on convergence for most nonlinear Kalman-type filters it is difficult to analyze the convergence of our method.  At each step of the algorithm we estimate the local average of the observation model error from the previous estimates $\hat b^{(\ell)}_k$ and then add this estimate to the observation function.  Notice that if the same state estimates $x_k^{(\ell+1)}=x_k^{(\ell)}$ were found in the next iteration of the Kalman filter, then the observation model error estimates would be unchanged.  Informally, if the state estimates only change by a small amount and if $g$ is continuous then the observation model error estimates should also only change by a relatively small amount.  In the next section we will present an interpretation of the method as an alternating minimization approach for estimating the local observation model error parameters.  Moreover, we will present numerical results demonstrating convergence for strongly nonlinear systems with extremely large error in the specification of the observation function.

\subsection{Interpretation as alternating minimization algorithm}

The method introduced above can be viewed as belonging to the family of projection algorithms in optimization theory called alternating minimization algorithms \cite{am1,am2}. Implicit to the above construction is the following nonparametric representation of the estimated global observation model error $b^{(\ell)}(x)$, which interpolates  the errors at each $x_k$ as
\[ b^{(\ell)}(x_k) =\sum_{i=1}^N \hat b_{k_i}^{(\ell)}  \frac{e^{-||z_{k_j}-z_k||/\sigma}}{\sum_{j=1}^{N}e^{-||z_{k_j}-z_k||/\sigma}} = \sum_{j=1}^N \hat b_{k_j}^{(\ell)} \frac{e^{-d(x,x_{k_j})/\sigma}}{\sum_{j=1}^N e^{-d(x,x_{k_j})/\sigma}}, \]
where $\{x_{k_j}\}_{j=1}^N$ are the $N$ nearest neighbors of the input $x$.  Takens' theorem \cite{takens,SYC} states that we can use the delay coordinate vectors $z_{k_j}$ as a proxy for the unknown true states $x_{k_j}$.  Using the Euclidean distance on the proxy vectors $z_{k_j}$ implicitly changes the distance function in state space to a metric $d$, which is consistent since all metric are equivalent in Euclidean space, and this has really only affected the weights in the average.  Notice that the finite set of parameters $\{\hat b_k^{(\ell)}\}$ determine the function $b^{(\ell)}(x)$.  From \eqref{e2} we assume that
\[ y_k = g(x_k) + b(x_k) + \nu_k \]
where $\nu_k$ is mean zero Gaussian noise with covariance matrix $R$.  Thus, the likelihood of the estimated observation model error $ b^{(\ell)}(x)$ can be estimated on the data set as
\begin{equation}\label{toMax} P\left(x_k^{(\ell)} \, | \,  b^{(\ell)}\right) \propto \prod_{k=1}^T \exp\left(-\frac{1}{2}||y_k-g(x_k^{(\ell)})- b^{(\ell)}(x_k^{(\ell)})||_{R}^2-\frac{1}{2}||x_{k+1}^{(\ell)}-f(x_k^{(\ell)})||_{Q}^2\right) \end{equation}
where $||\nu||_{R}^2 = \nu^\top R^{-1} \nu$ is the norm induced by the covariance matrix $R$.  Our goal is to maximize the probability simultaneously with respect to both the state estimate $x_k^{(\ell)}$ and the observation model error estimate $\hat b^{(\ell)}$, or equivalently, to minimize $-\log P$, the negative log likelihood. 

At the $\ell$-th step of our approach, we first fix the observation model error estimate $ b^{(\ell)}$ and use the nonlinear Kalman filter to approximate the best estimate of the state $x_k^{(\ell)}$ given the current estimate of the observation model error.  The nonlinear Kalman filter is approximating the solution which maximizes \eqref{toMax} where $b^{(\ell)}$ is fixed.  One could also apply a variational filtering method to achieve this maximization.

\begin{figure}[h]
\begin{center}
\subfigure[]{\includegraphics[width=0.32\columnwidth]{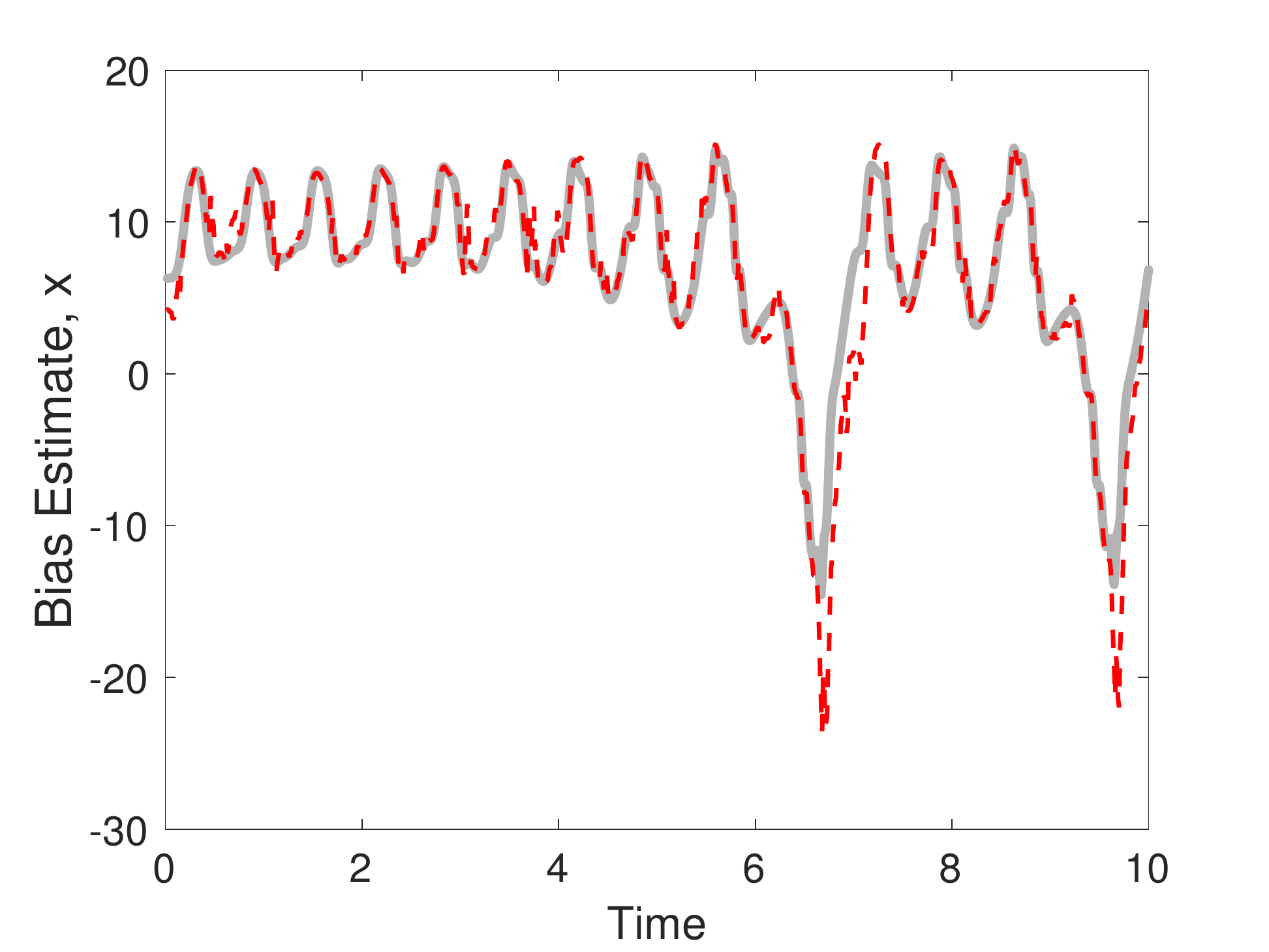}}
\subfigure[]{\includegraphics[width=0.32\columnwidth]{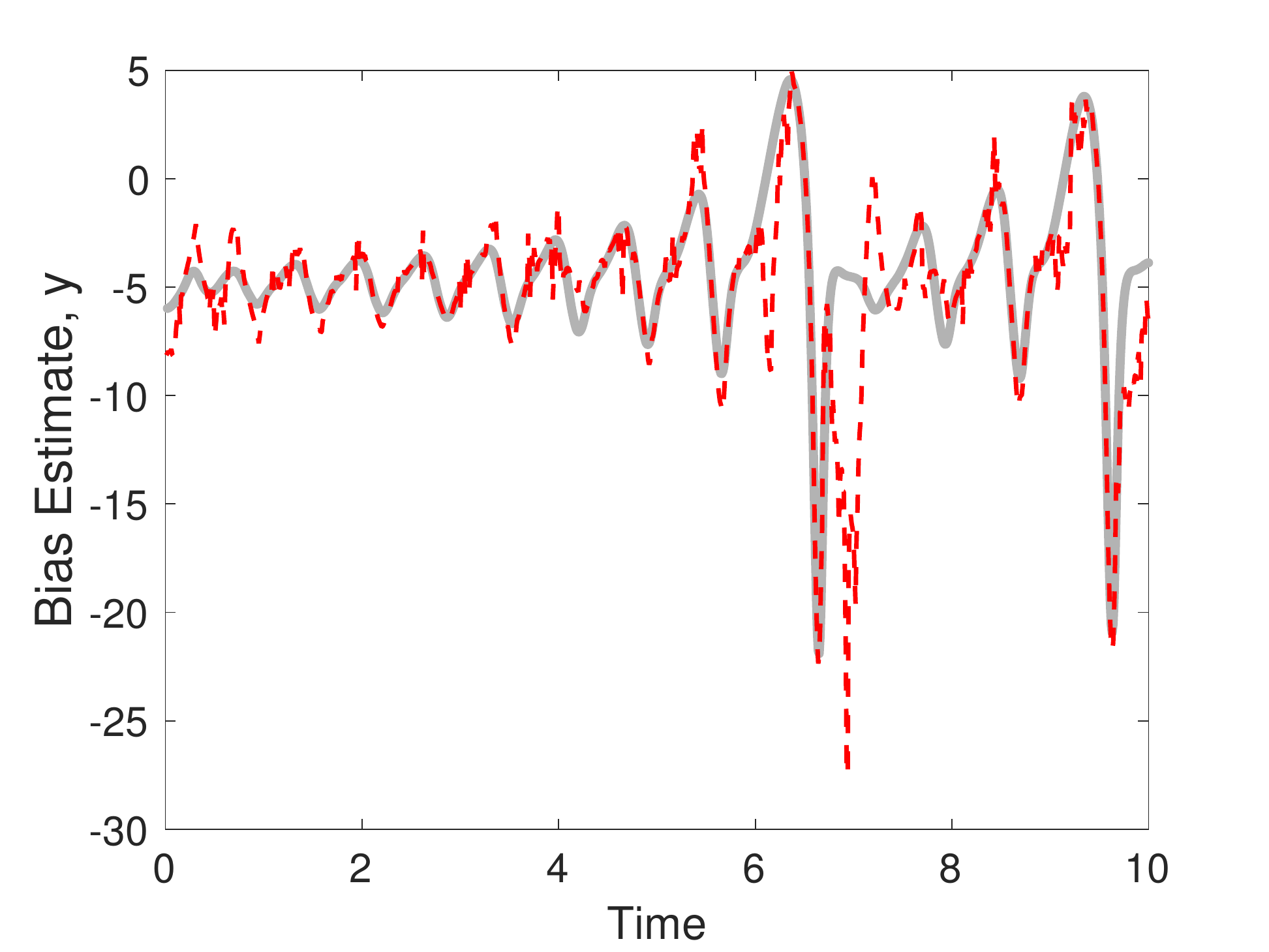}}
\subfigure[]{\includegraphics[width=0.32\columnwidth]{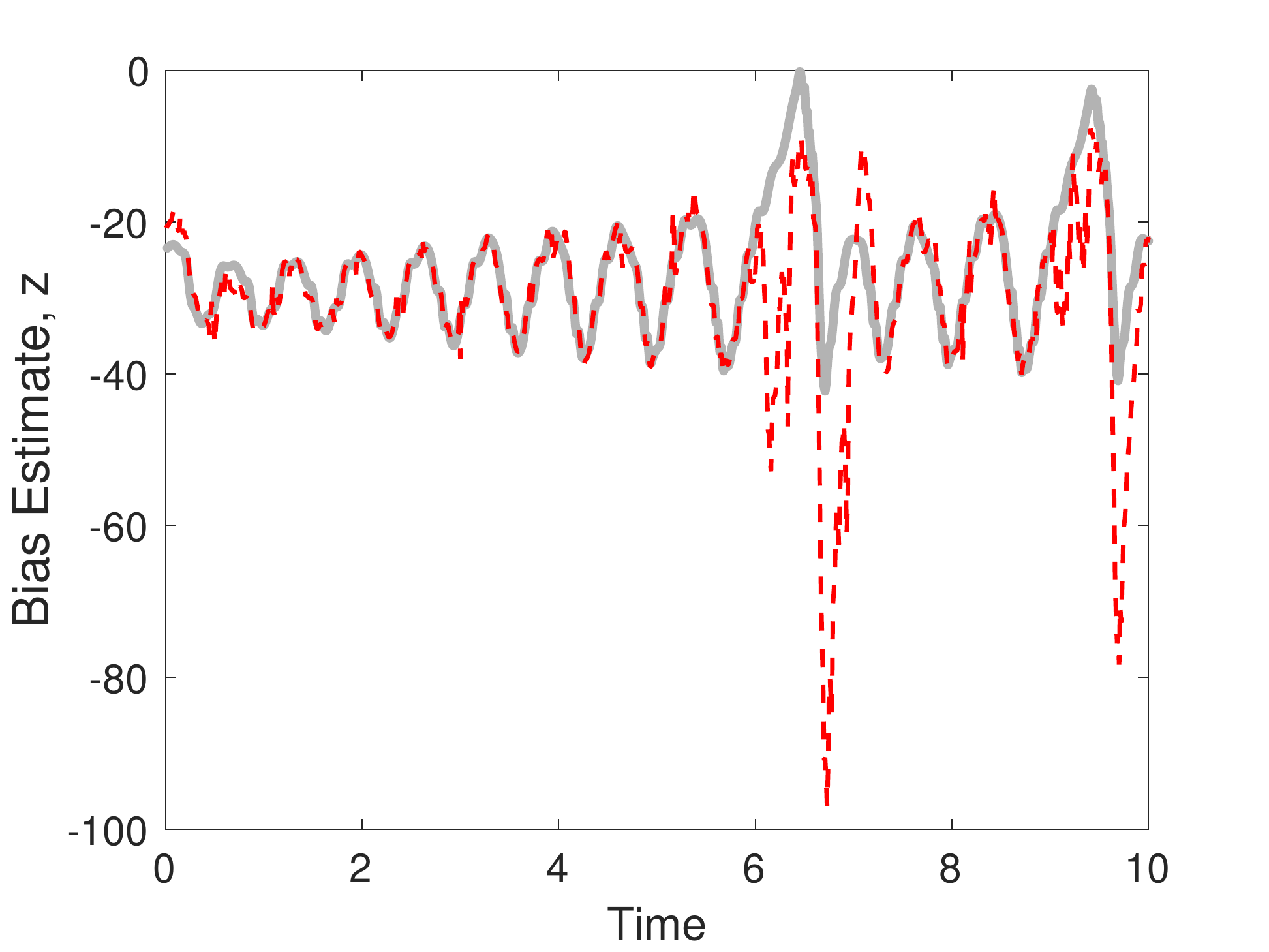}}
\end{center}
\caption{Comparison of the observation model error correction which solves \eqref{biasReg} (red, dashed) to the correction given by \eqref{biasSimple} (grey, solid) which is used in all the examples below.  Observation errors are shown from the Lorenz-63 example described below (see Fig.~\ref{figure3}).}
\label{BiasComp}
\end{figure}

Next, we fix the estimate $x_k^{(\ell)}$ and estimate the parameters $\hat b_k^{(\ell+1)}$ to maximize \eqref{toMax}.  Since the second term in the exponential is independent of $\hat b_k^{(\ell+1)}$, the solution which maximizes \eqref{toMax} is simply the solution to the linear system of equations
\begin{equation}\label{biasReg} y_k-g(x_k^{(\ell)}) =  b^{(\ell)}(x_k^{(\ell)}) =\sum_{j=1}^N \hat b_{k_j}^{(\ell)} \frac{e^{-d(x,x_{k_j})/\sigma}}{\sum_{j=1}^N e^{-d(x,x_{k_j})/\sigma}}.\end{equation}
Instead of explicitly solving this system, in our implementation we simply used the approximate solution given by 
\begin{equation}\label{biasSimple} \hat b_k^{(\ell)} = y_k-g(x_k^{(\ell)}) \end{equation}
since each point is its own nearest neighbor and $d_{k_1}=0$ yields the largest weight in the summation.  In Fig.~\ref{BiasComp} we show that the observation model error estimates \eqref{biasReg} and \eqref{biasSimple} are very similar, but \eqref{biasSimple} is much faster to compute and is more numerically stable so we will use \eqref{biasSimple} in all the examples below.

\subsection{Ensemble Kalman filtering with observation model error correction}
In this section we assume a nonlinear system with $n$-dimensional state vector $x$ and $m$-dimensional observation vector $y$ defined by \eqref{e1} as described above.
The ensemble Kalman filter (EnKF) approximates a nonlinear system by forming an ensemble, such as through the unscented transformation (see for example \cite{simon}). At the $k$th step of the filter there is an estimate of the state $x^+_{k-1}$ and the covariance matrix $P^+_{k-1}$. In the unscented version of the EnKF, the singular value decomposition is used to find the symmetric positive definite square root $S^+_{k-1}$ of the matrix $P^+_{k-1}$, allowing us to form an ensemble of $E$ state vectors where the $i^{th}$ ensemble member is denoted $x_{i,k-1}^+$.

The EnKF alternates between a forecast, where the state is predicted, and an analysis, where the current observation is used to correct the state prediction. The model $f$ is applied to the ensemble, advancing it forward one time step, and then observed with function $g^{(\ell)}$
\begin{eqnarray}
\label{e2a}
x_{i,k}^- &=& f\left(x_{i,k-1}^+\right) \nonumber\\
y_{i,k}^- &=& g^{(\ell)}\left(x_{i,k}^-\right).
\end{eqnarray}
Notice that in the ideal filtering situation we would apply the true observation function $h$ in \eqref{e2a}, however, in this context we assume that we are only given an incorrect observation function $g$.  In the initial iteration of the filter ($\ell=0$) we simply use the best available observation function $g^{(0)}=g$, and in future iterations $(\ell>0)$ we incorporate the $\ell$-th observation model error estimate to form $g^{(\ell)} = g + \hat b^{(\ell)}$ as described above.  Notice that each ensemble member has the same correction $\hat b^{(\ell)}$ applied since the correction is computed based on the neighbors in delay-embedded observation space, so the neighbors do not change based on the state estimate or iteration of the algorithm.  We emphasize that the state estimate and observation model error estimates change at each iteration, but the indices of the neighbors, $k_1,...,k_N$ that are used to estimate the observation model error at time step $k$ do not change (they are independent of $\ell$).

The mean of the resulting state ensemble gives the prior state estimate $x_k^-$ and the mean of the observed ensemble is the predicted observation $y_k^-$.  Denoting the covariance matrices $P_k^-$ and $P_k^y$ of the resulting state and observed ensemble, and the cross-covariance matrix $P_k^{xy}$ between the state and observed ensembles, we define
\begin{eqnarray}
P_k^- &=& \frac{1}{E}\sum_{i= 1}^{E} \left(x_{i,k}^-  -x_k^-\right) \left(x_{i,k}^- -x_k^-\right)^T + Q\nonumber\\
P_k^y &=&  \frac{1}{E}\sum_{i= 1}^{E} \left(y_{i,k}^-  -y_k^-  \right) \left(y_{i,k}^-  -y_k^- \right)^T + R \nonumber\\
P_k^{xy} &=& \frac{1}{E}\sum_{i= 1}^{E} \left(x_{i,k}^-  -x_k^- \right) \left(y_{i,k}^-  -y_k^-  \right) ^T
\end{eqnarray}
and use the equations
\begin{eqnarray} 
K_k &=& P^{xy}_k(P^{y}_k)^{-1}\nonumber\\
P^{+}_k &=& P^{-}_k-K_kP^{yx}_k\nonumber\\
{x}^{+}_k &=& {x}^{-}_k+K_k\left(y_k-{y}_k^- \right).
\end{eqnarray}
to update the state and covariance estimates with the observation $y_k$. $Q$ and $R$ are generally unknown quantities that have to be estimated, an area known as adaptive filtering. 

In this article, we use the method of \cite{berry2} for the adaptive estimation of these noise covariance matrices.  This is a key component in our method since the $R$ covariance will be inflated by the adaptive filter to represent the error between the true observation function $h$ and the observation function $g^{(\ell)}$ that we actually use in the filter.  In other words, the adaptive filter is combining the covariance of the observation model error and the instrument noise into the $R$ covariance matrix.  As we iterate the algorithm (as $\ell$ increases) we find that $g^{(\ell)}$ more closely approximates the true observation function $h$ and the adaptive filter will find smaller values for $R$.

\section{Assimilating Lorenz-63 with an incorrect observation model}\label{L63section}
In the results presented below, we assume noisy observations are available from a system of interest and we implement an ensemble Kalman filter (EnKF) for state estimation.  The EnKF approximates a nonlinear system by forming an ensemble, such as through the unscented transformation (see for example \cite{simon}). Additionally, we use the method of \cite{berry2} for the adaptive estimation of the filter noise covariance matrices $Q$ and $R$. The correct observation function $h$ that maps the state to observation space is unknown, and in its place an incorrect function $g$ is chosen for use by the EnKF. Throughout, we will compare our corrected filter with the standard filter (essentially, the $\ell = 0$ iteration) which assumes no correction.


\begin{figure}[h]
\begin{center}
\subfigure[]{\includegraphics[width=0.4\columnwidth]{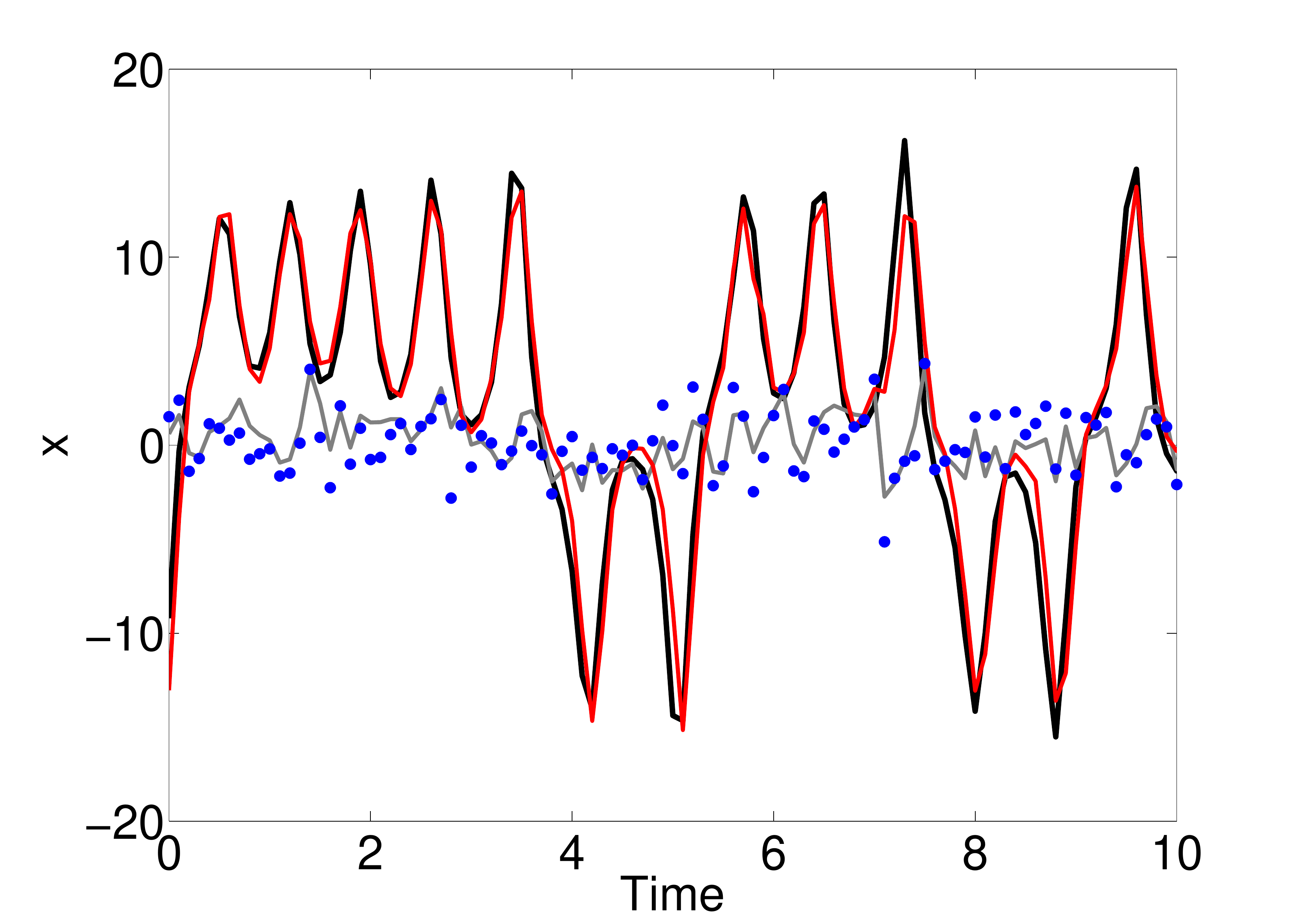}}
\subfigure[]{\includegraphics[width=0.4\columnwidth]{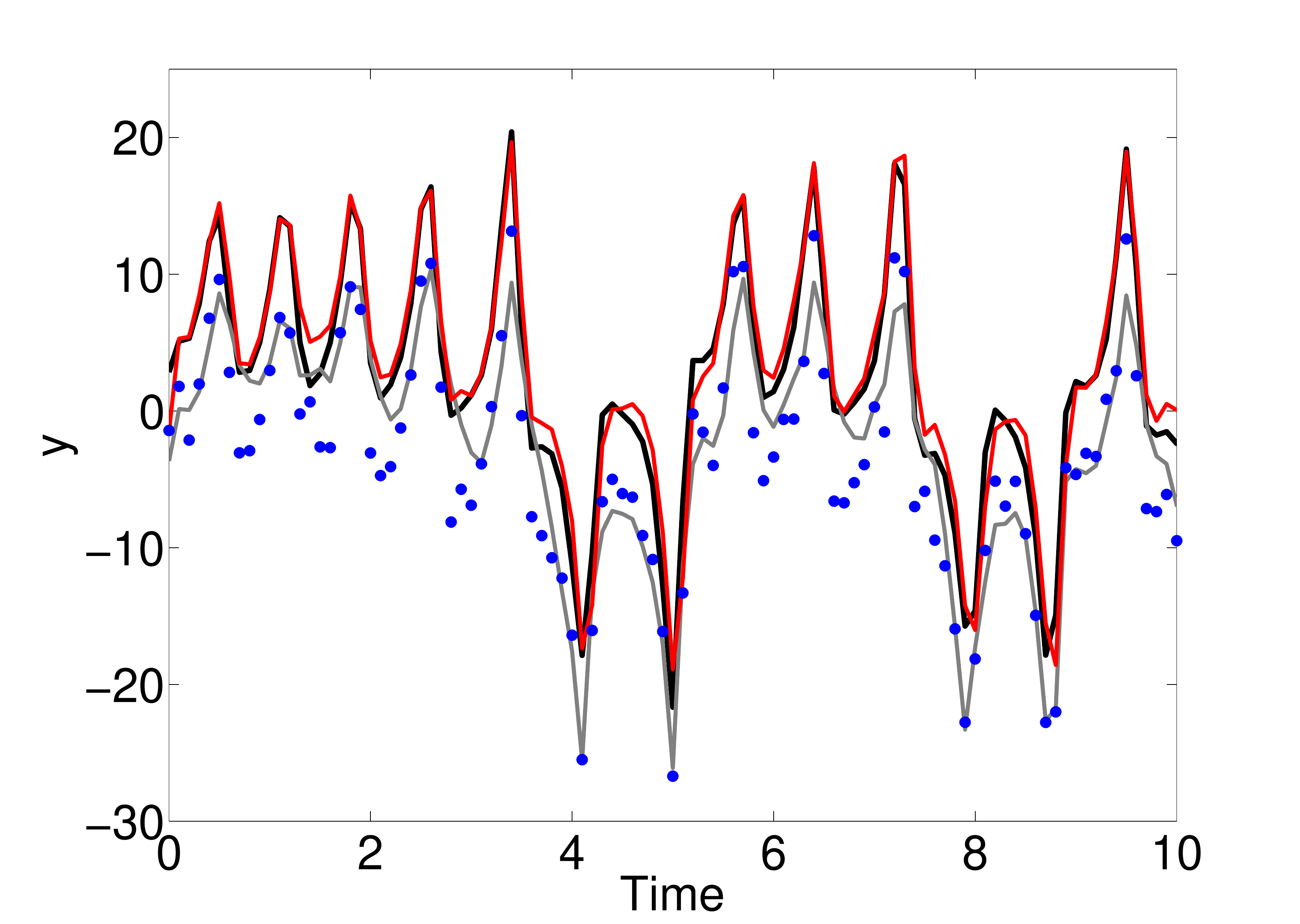}}\\
\subfigure[]{\includegraphics[width=0.4\columnwidth]{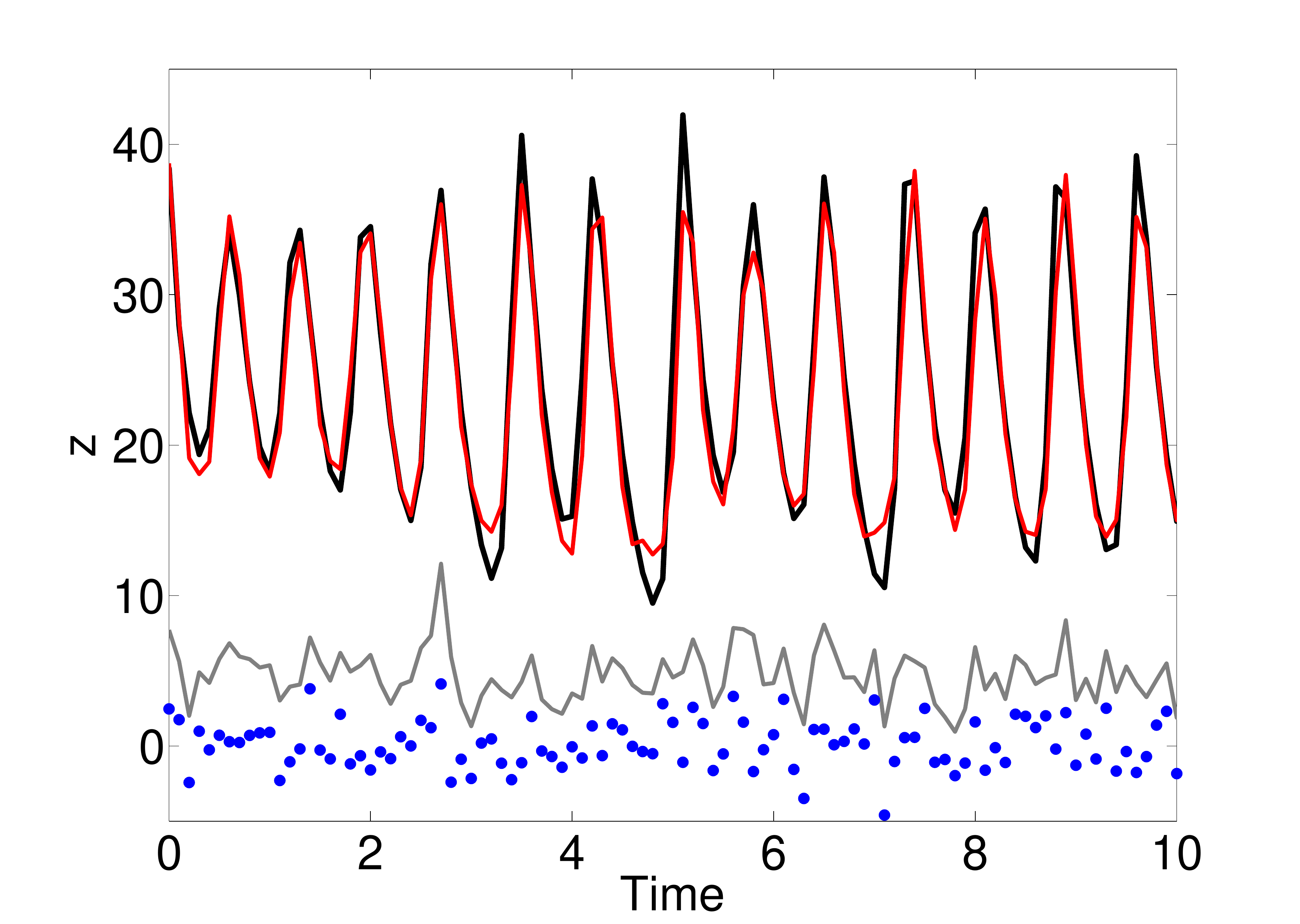}}
\subfigure[]{\includegraphics[width=0.4\columnwidth]{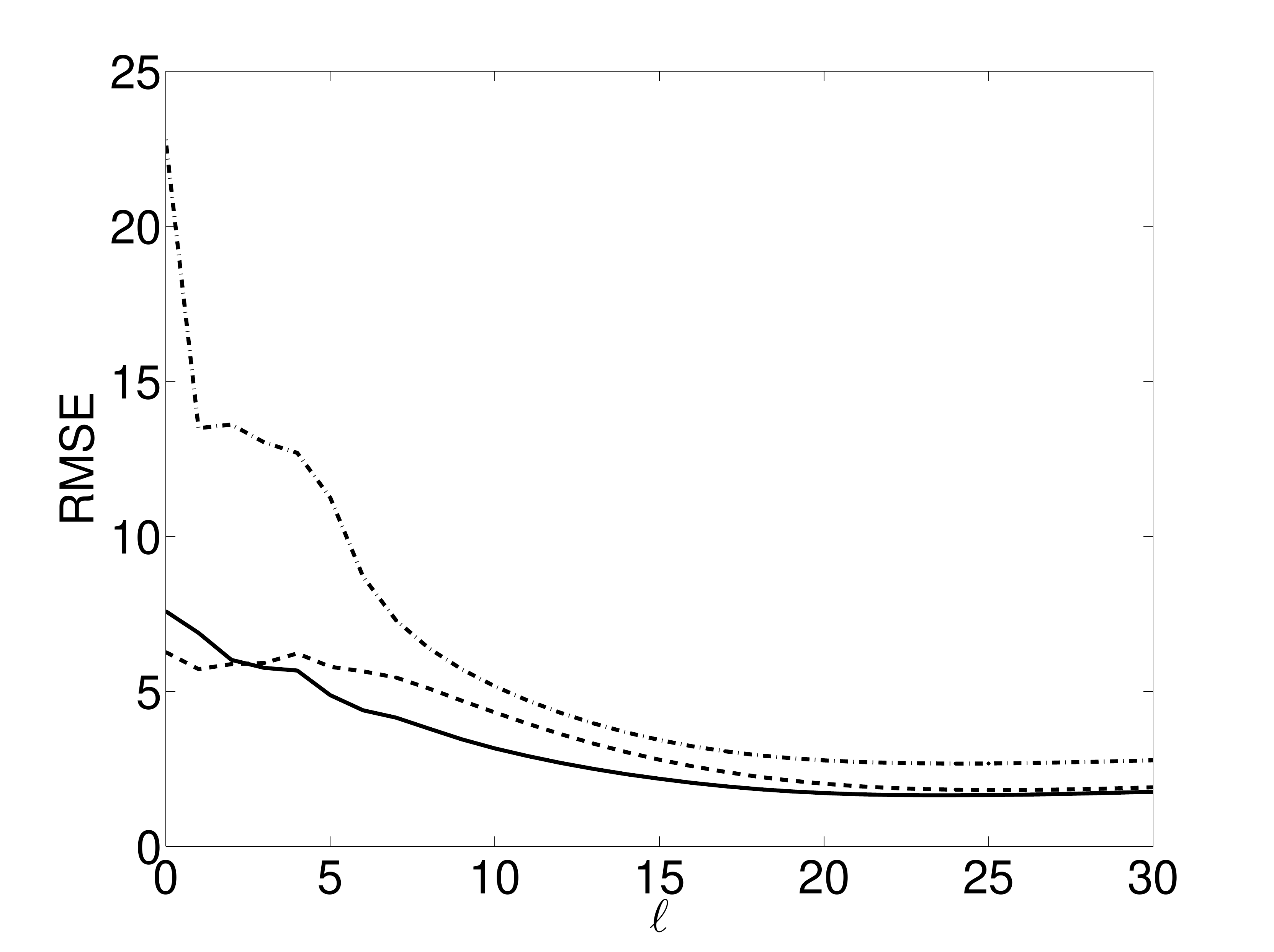}}
\end{center}
\caption{Results of filtering noisy Lorenz-63 (a) $x_1$ (b) $x_2$ and (c) $x_3$ time series when true observation function, $h$, is unknown and $R=2I_{3\times 3}$.  Notice the large difference between the true observations $h(\vec x_k)+\nu_k$ (blue circles) to the true state variables (solid black curve).  We compare the EnKF estimate using the wrong observation function, $g$, without observation model error correction (solid gray lines) and the EnKF estimate with correction (solid red lines) shown. 
(d) Plot of RMSE vs. iteration of the observation model error correction method, where $\ell = 0$ corresponds to the standard EnKF without correction. RMSE for $x$ (solid black line), $y$ (dashed black line) and $z$ (dotted black line) shown. After a sufficient number of iterations, the observation model error estimates converge as does the RMSE of the state estimate.}
\label{figure3}
\end{figure}

As a demonstrative example we consider the Lorenz-63 system \cite{Lorenz63}
\begin{eqnarray} \label{e6}
\dot{x_1} &=& \sigma(x_2-x_1)\nonumber\\
\dot{x_2} &=& x_1(\rho-x_3)-x_2 \\
\dot{x_3} &=& x_1 x_2-\beta x_3 \nonumber
\end{eqnarray}
where $\sigma = 10$, $\rho = 28$, $\beta = 8/3$.  We will assimilate 8000 noisy observations of the system, sampled at rate $dt = 0.1$, and corrupted by independent Gaussian observational noise, $\nu_k$, with mean zero and covariance $R=2I_{3\times 3}$.  Our goal is to filter the observations
\[ \vec y = h(\vec x) + \nu_k \]
(see Fig.~\ref{figure3}, blue circles) and reconstruct the underlying state, $\vec x$, (Fig.~\ref{figure3}, solid black lines). However, we assume that the true observation function $h$, given by 
\begin{eqnarray*}
h(\vec x) = h\left(\left[\begin{array}{c} x_1\\ x_2 \\ x_3 \end{array} \right] \right) = \left[\begin{array}{c} \sin(x_1)\\x_2-6\\ \cos(x_3)\end{array}\right]
\end{eqnarray*}
is unknown to us.  Instead, the EnKF will use an incorrect observation function $g$, given by
\begin{eqnarray*}
g(\vec x) = g\left(\left[\begin{array}{c} x_1\\ x_2 \\ x_3 \end{array} \right] \right) = \left[\begin{array}{c} x_1\\ x_2 \\ x_3\end{array}\right].
\end{eqnarray*}
Using the incorrect mapping $g$, and with no estimate of the observation model error, the filter's reconstruction of the system state suffers substantially (Fig.~\ref{figure3}(a)-(c), solid gray lines).  We should note that even obtaining these poor estimates requires adaptive estimation of the system and observation noise covariance matrices $Q$ and $R$ used by the EnKF.  The RMSE for reconstructing the three Lorenz-63 variables $x_1,x_2$ and $x_3$ using an EnKF with observation function $g$ and no observation model error correction is 8.10, 6.77 and 22.33 respectively. This is not surprising, since without the correct observation function the analysis step of the EnKF, where the state and covariance estimates are updated, suffers due to the errors in mapping the predicted state to observation space.

Using our proposed method, the EnKF state estimate can be improved by iteratively building an approximation of the observation model error, essentially augmenting our observation function. In building our reconstruction of the observation model error, we use $d =  2$ delays and $N = 100$ nearest neighbors. After $M = 20$ iterations of our method, we are able to obtain and accurate estimate of the Lorenz-63 state (Fig.~\ref{figure3}(a)-(c), solid red lines). The resulting error in our estimates is significantly smaller (RMSE of 2.11, 1.77 and 2.91 for $x$, $y$ and $z$ respectively) compared to filtering without an observation model error correction.

Fig.~\ref{figure3}(d) shows the error in our estimation of $x$ (solid black line), $y$ (dashed black line) and $z$ (dotted black line) as a function of number of iterations of our algorithm. We note that $\ell = 0$ corresponds to running the EnKF without any observation model error. At each iteration, we obtain a better reconstruction of the observation model error which helps improve our estimate of the state in the next iteration. At a certain point, our reconstruction of the observation model error and system state converges, a period indicated by the plateau in our RMSE plot.

\begin{figure} [ht]
\center
\subfigure[]{\includegraphics[width=0.45\columnwidth]{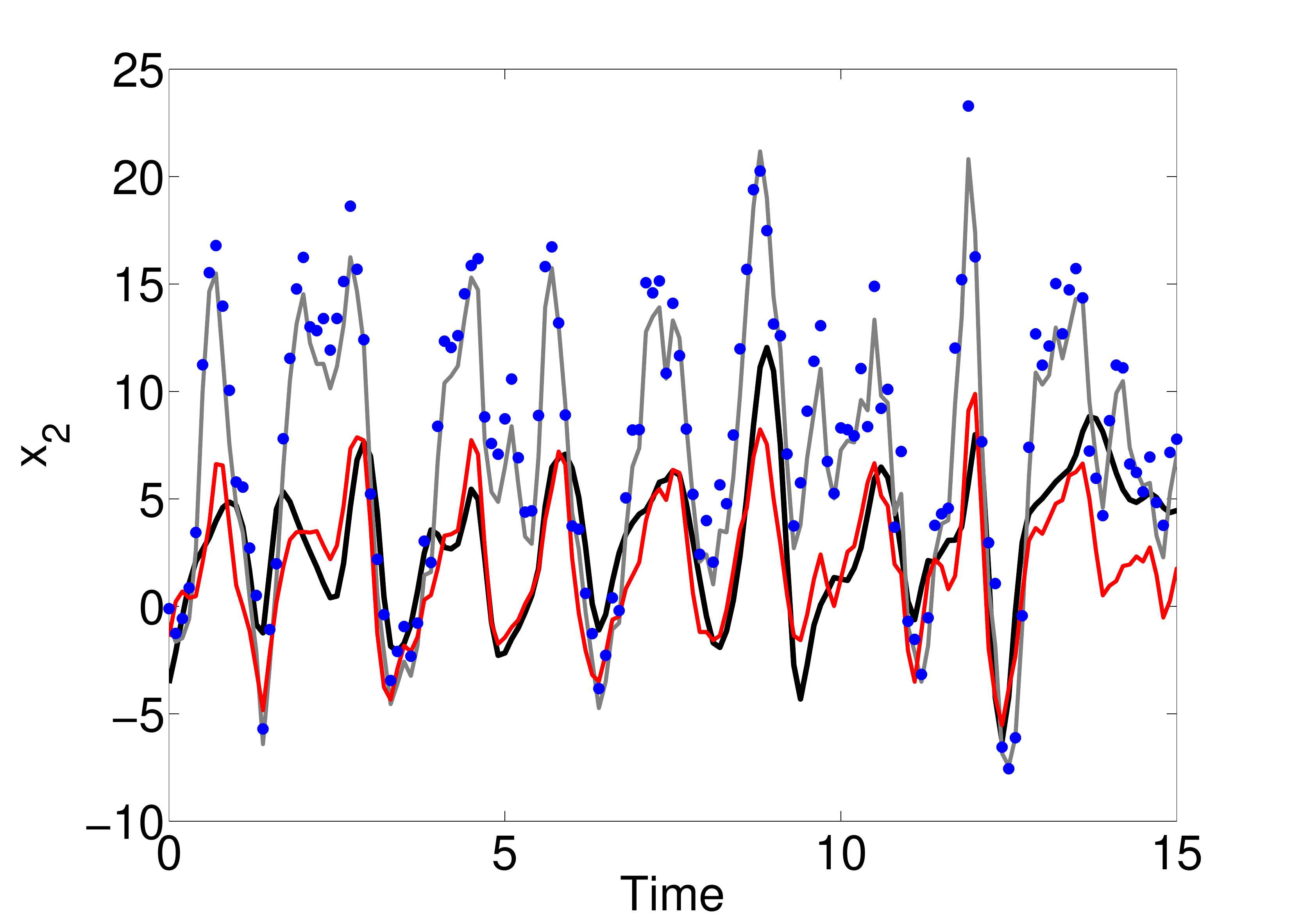}}
\subfigure[]{\includegraphics[width=0.45\columnwidth]{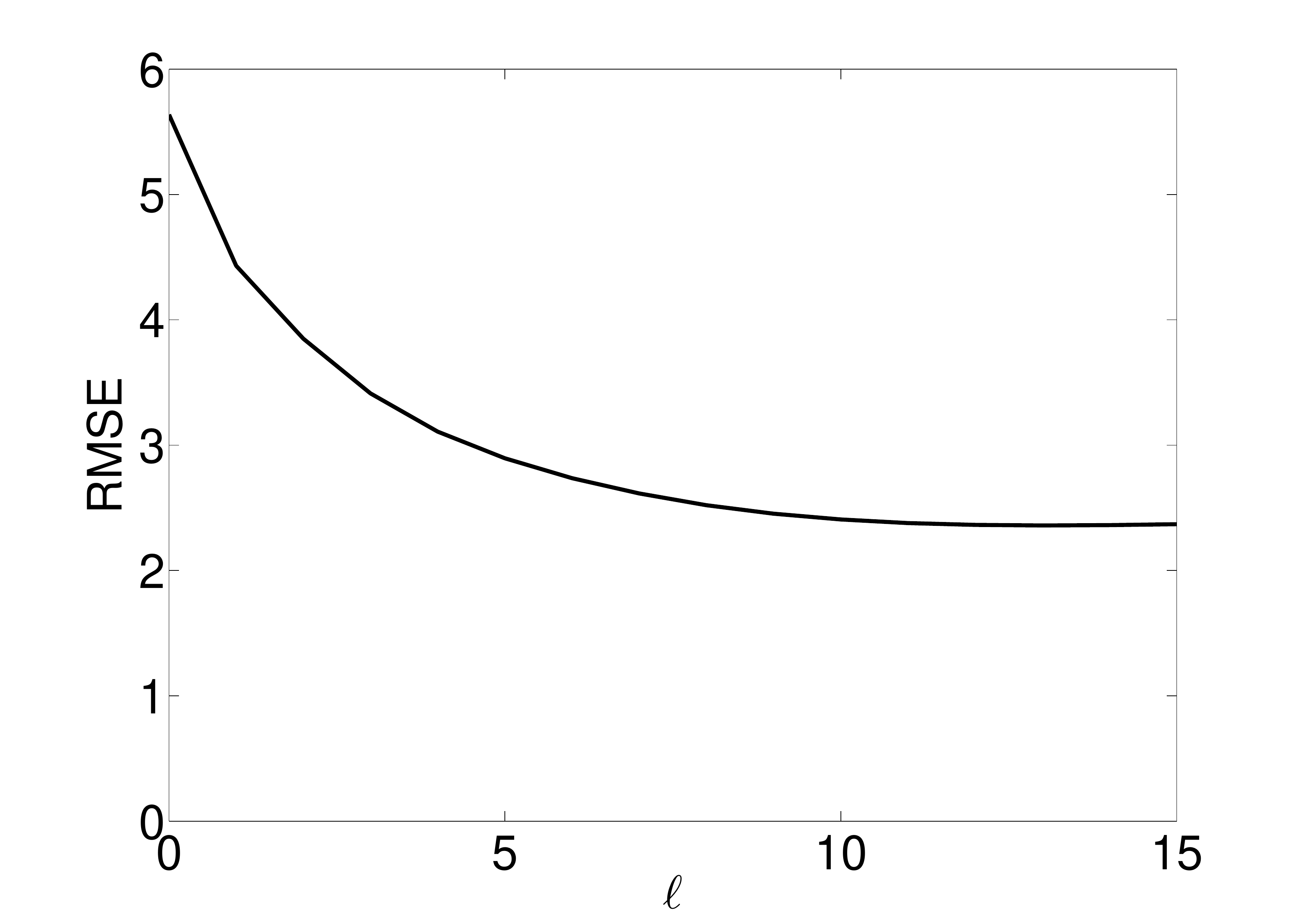}}
\caption{Results of filtering a noisy 10 dimensional Lorenz- 96 ring when the true observation function is unknown. (a) Representative results demonstrated by the $x_{2}$ node. We filter the noisy observation (blue circles) in an attempt to reconstruct the underling state (solid black line). Without observation model error correction, the EnKF estimate (solid gray line) is unable to track the true state (RMSE = 5.83). With observation model error correction (solid red line), our estimate of the state improves substantially (RMSE = 2.37). (b) Average RMSE of Lorenz-96 ring as a function of iteration shown. Similarly to the previous example, after a sufficient number of iterations our method converges to an estimate of the observation model error and system state, demonstrated by the convergence of the RMSE.}
\label{figure4}
\end{figure}

\section{Spatiotemporal observation model error correction}

We now consider a coupled ring of $K$ nodes of Lorenz-96 \cite{lorenz1996predictability} equations
\begin{eqnarray} \label{e7}
\dot{x}_i = (ax_{i+1}-x_{i-2})x_{i-1}-x_i+F
\end{eqnarray}
where $a = 1$ and $F = 8$. The Lorenz-96 system is a convenient example since it allows for a range of higher dimensional complex behavior by adjusting the number of nodes in the system. We assume that 10000 observations, corrupted by mean 0 Gaussian noise with variance equal to 2, are available from each node in the ring. Denoting $\mathbf{x} =[ x_1,x_2,\hdots,x_{K}]$, the true observation function $h$ for this system is defined as $h\left(\mathbf{x}\right) = C\mathbf{x}$, where
\begin{eqnarray*}
C = \left[\begin{array}{cccccccccccc} c_1 & c_2 & 0 & \cdots & \cdots & \cdots & \cdots & c_3\\
c_3 & c_1 & c_2 & 0 & & & & \vdots\\
0 & c_3 & c_1 & c_2 & \ddots & & & \vdots\\
\vdots & 0 & \ddots & \ddots & \ddots & \ddots & & \vdots\\
\vdots & & \ddots & \ddots & \ddots & \ddots & 0 & \vdots\\
\vdots & & & \ddots & c_3 & c_1 & c_2 & 0\\
\vdots & & & & 0 & c_3 & c_1 & c_2\\
c_2 & \cdots & \cdots  & \cdots & \cdots & 0 & c_3 & c_1\end{array}\right],
\end{eqnarray*}
 $c_1 = 1,c_2 = 1.2,c_3 = 1.1$. In effect, our observations at each node in the ring is a linear combination of the current node and its two spatial neighbors. The true observational mapping $h$ is unknown to us, and in its place we assume the incorrect function
\begin{eqnarray*}
g\left(\mathbf{x}\right) = \mathbf{I}_{K\times K}\mathbf{x}
\end{eqnarray*}
where $I_{K\times K}$ is the $K \times K$ identity matrix.

We first consider a $K= 10$ dimensional Lorenz-96 ring. Fig.~\ref{figure4} shows the results of reconstructing the 10 dimensional Lorenz-96 state. Fig.~\ref{figure4}(a) shows a representative reconstruction of the $x_{2}$ state (similar results are obtained for each node of the ring). Given the noisy observations (blue circles), the EnKF without observation model error correction (solid gray line) is unable to estimate the true trajectory (solid black line), resulting in an RMSE of 5.83. Accounting for the observation model error ($M = 15$ iterations, $d = 2$ delays and $N = 100$ neighbors), we are able to improve our reconstruction of the $x_{2}$ trajectory (solid red line, RMSE = 2.37). Similarly as in the Lorenz-63 example, we see in Fig.~\ref{figure4}(b) that as the number of iterations of our observation model error correction method increases we eventually converge to a stable RMSE.

\begin{figure}[ht]
\center
\includegraphics[width=0.8\columnwidth]{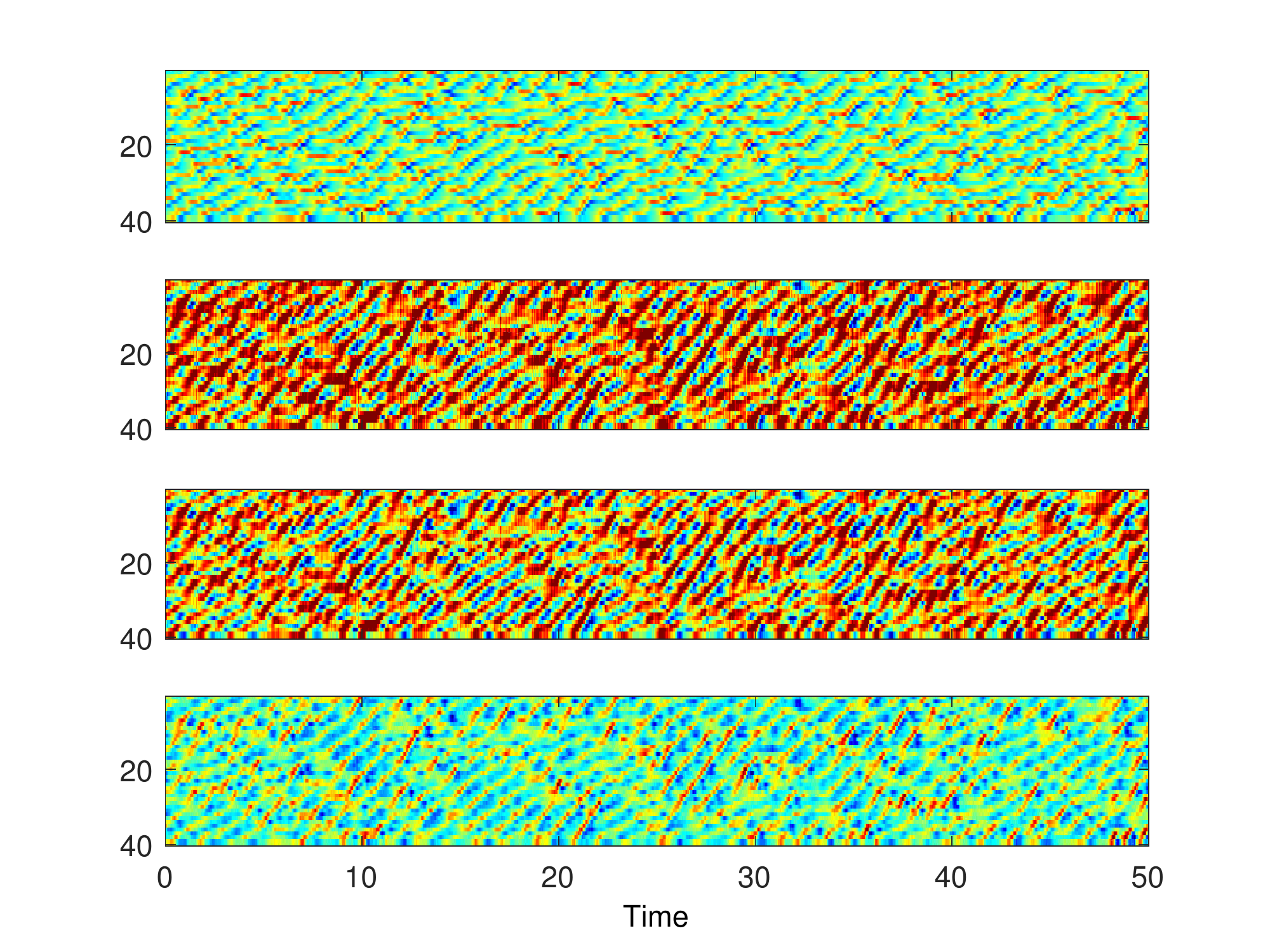}
\caption{Results of filtering a noisy 40 dimensional Lorenz-96 system. True spatiotemporal dynamics (top), noisy observations of the system (second plot), estimate without observation model error correction (third plot) and estimate with observation model error correction (bottom plot) shown. Without correction, we obtain a poor estimate of the system dynamics (average RMSE = 5.12). With correction, our estimate is improved (average RMSE = 2.50).}
\label{figure5}
\end{figure}

Given the success on the smaller Lorenz-96 system, we now consider a $K = 40$ dimensional ring. Fig.~\ref{figure5} shows the spatiotemporal plots of the system. The top plot shows the true system dynamics and the second plot our noisy observations of the system. Similarly to the 10 dimensional ring, the filtering without observation model error correction is unable to provide an accurate reconstruction of the system state (third plot). The high dimensionality of the system can make finding accurate nearest neighbors for observation model error reconstruction difficult. We implement a spatial localization technique when finding neighbors, whereby for each node we look for neighbors in a delay-coordinate space consisting of its delays and the delays of its two spatial neighbors. While our method can be successfully implemented in this high dimensional example without localization, results are improved through use of the localization technique. The bottom plot of Fig.~\ref{figure5} shows the resulting filter estimate with observation model error correction. Again, we see that there is a substantial improvement in the state reconstruction and we are able to obtain a more accurate representation of the true system dynamics.

\section{Correcting error in cloudy satellite-like observations without training data}

The presence of clouds is a significant issue in assimilation of satellite observations. Clouds can introduce significant observation model error into the results of radiative transfer models (RTM).  As previously mentioned, a recently developed method \cite{berry2017} is able to learn a probabilistic observation model error correction using training data consisting of pairs of the true state and the corresponding observations.  Of course, requiring knowledge of the true state in the training data is a significant restriction, and while methods such as reanalysis or local large-scale data gathering are possible, it would be extremely advantageous to remove this requirement. The innovation of the method introduced here is that we do not require knowledge of the true state in the training data.  Instead, we use an iterative approach to learn local observation model error corrections based on delay reconstruction in observation space.  In this section we will apply our method to an RTM and show that the observation model error can be iteratively learned without the training data.

\begin{figure}[h]
\center
\subfigure[]{\includegraphics[width=0.24\columnwidth]{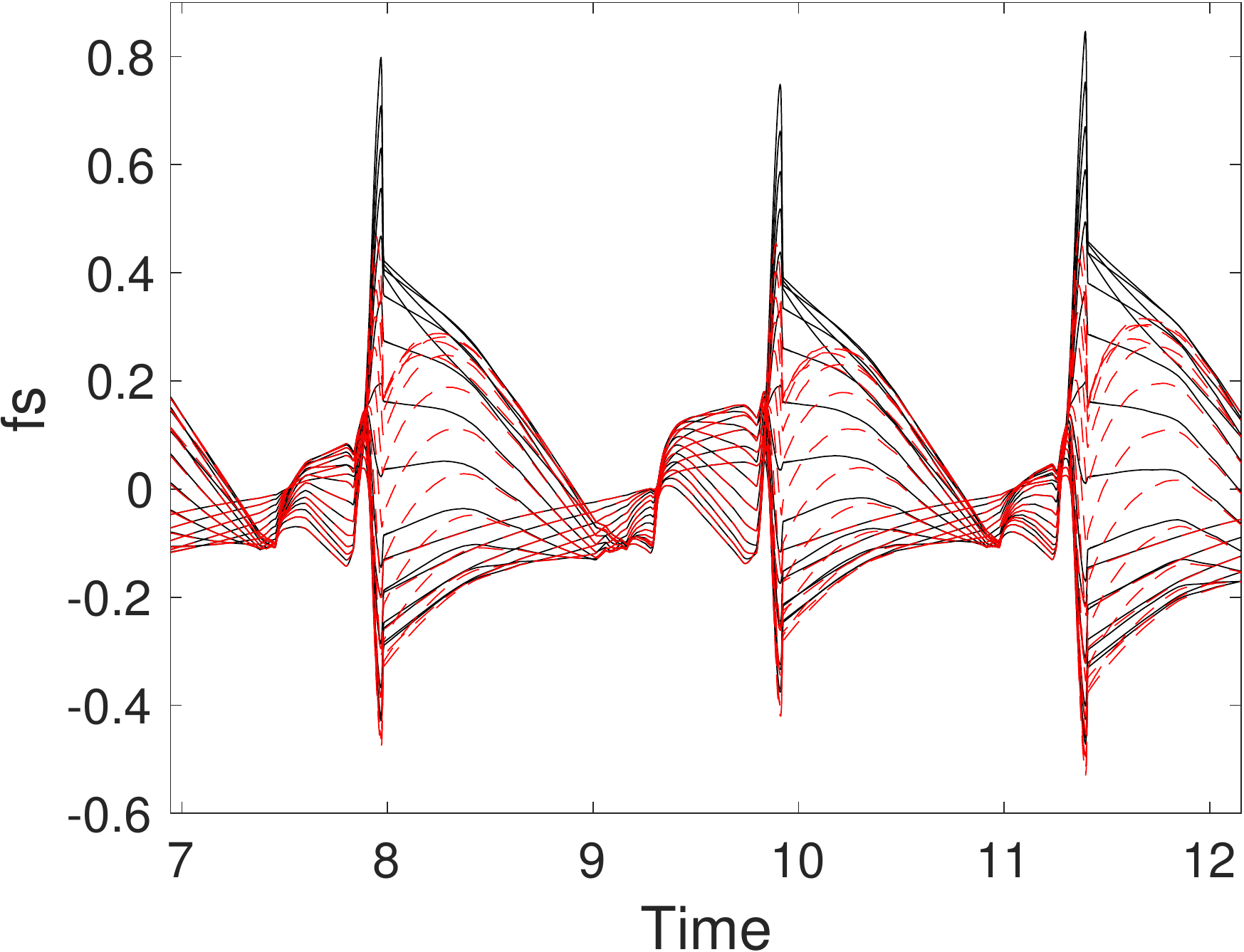}}
\subfigure[]{\includegraphics[width=0.24\columnwidth]{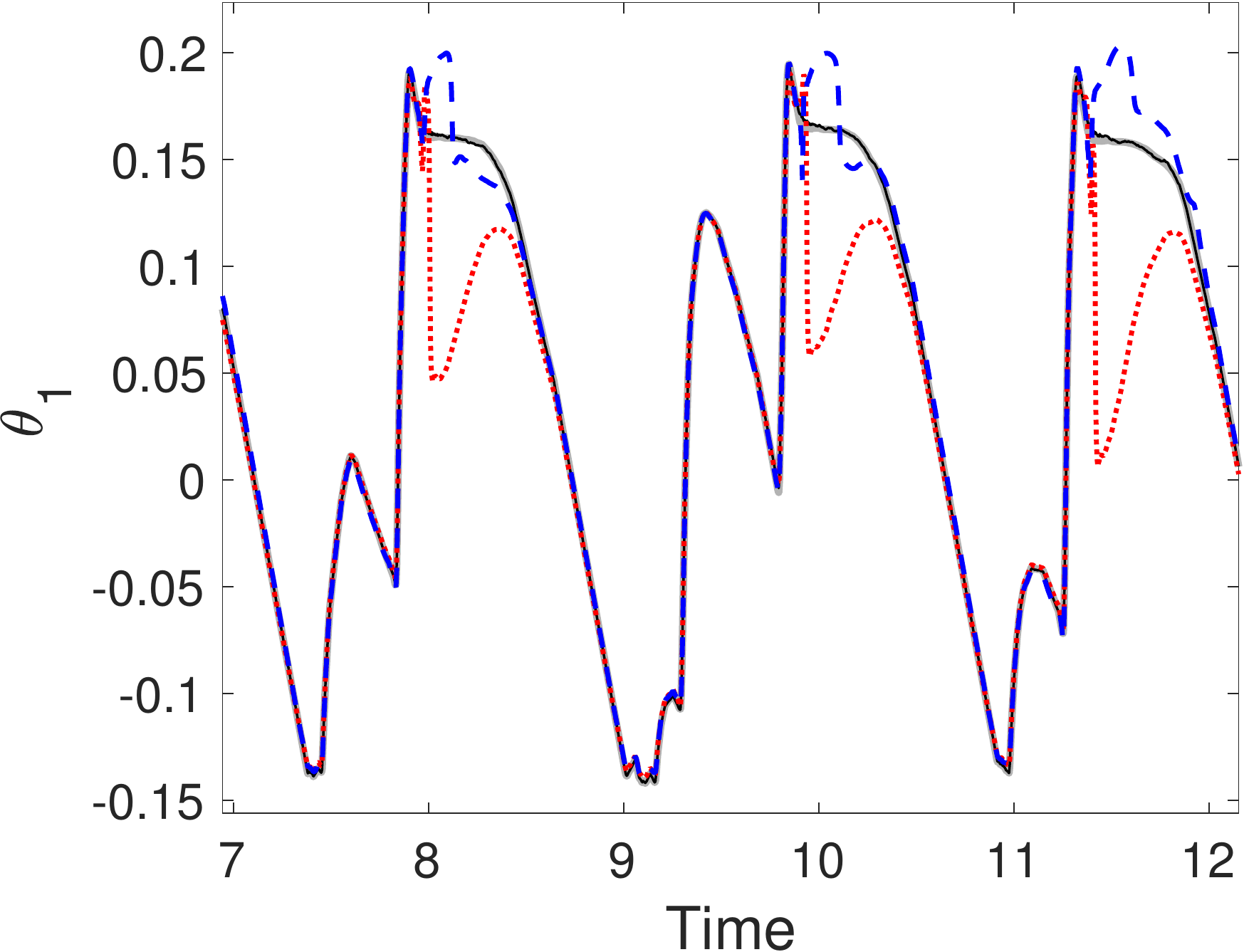}}
\subfigure[]{\includegraphics[width=0.24\columnwidth]{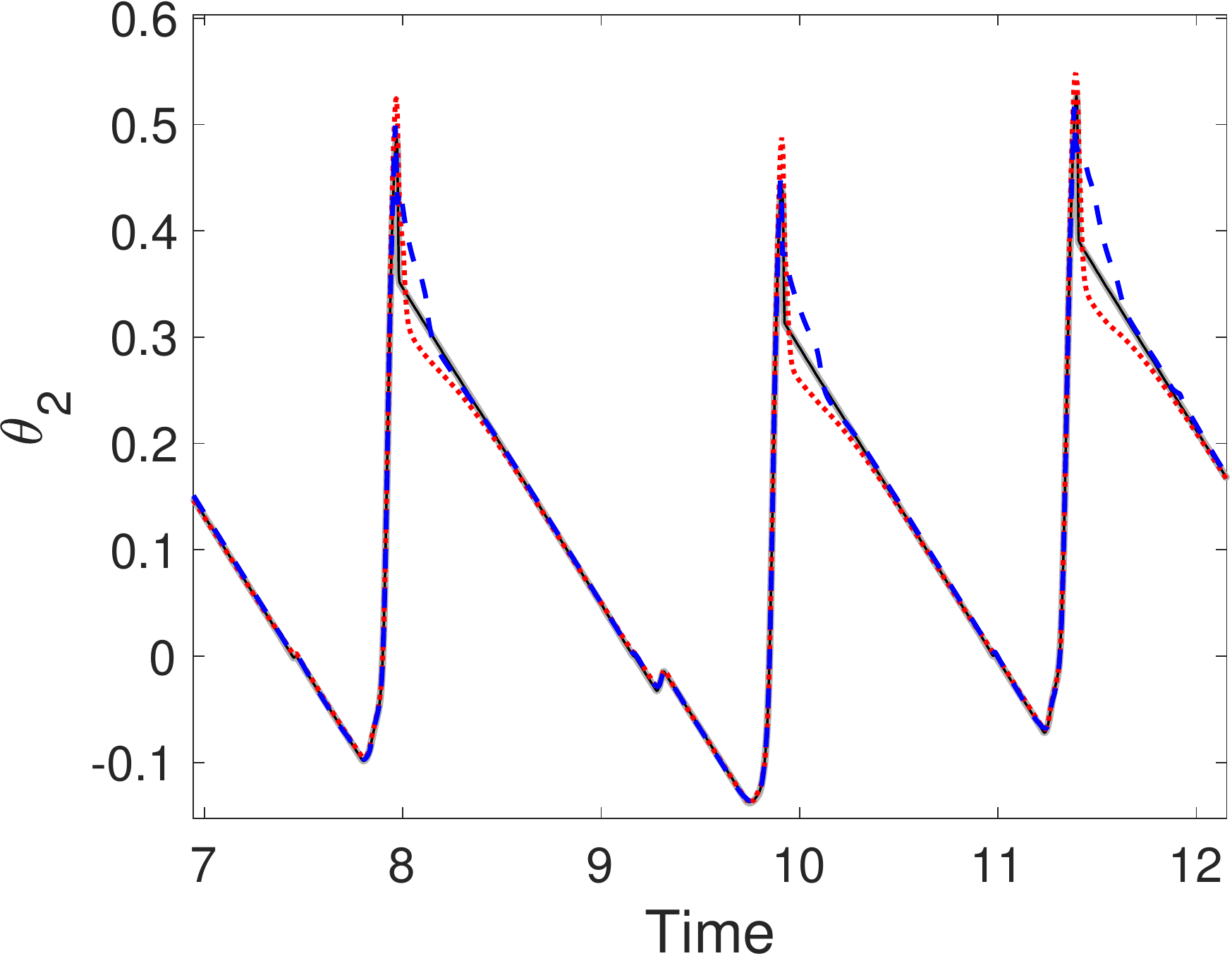}}
\subfigure[]{\includegraphics[width=0.24\columnwidth]{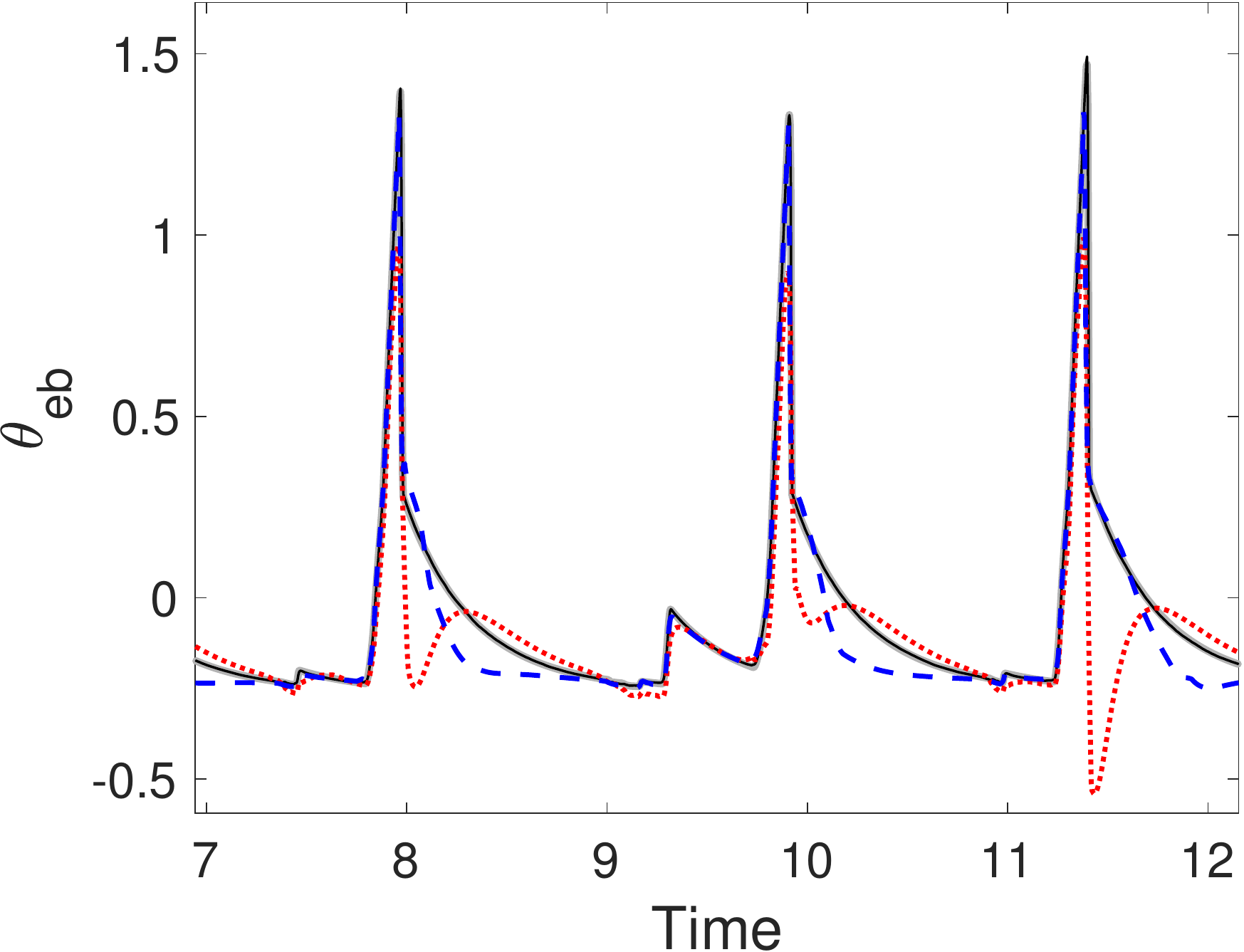}}
\subfigure[]{\includegraphics[width=0.24\columnwidth]{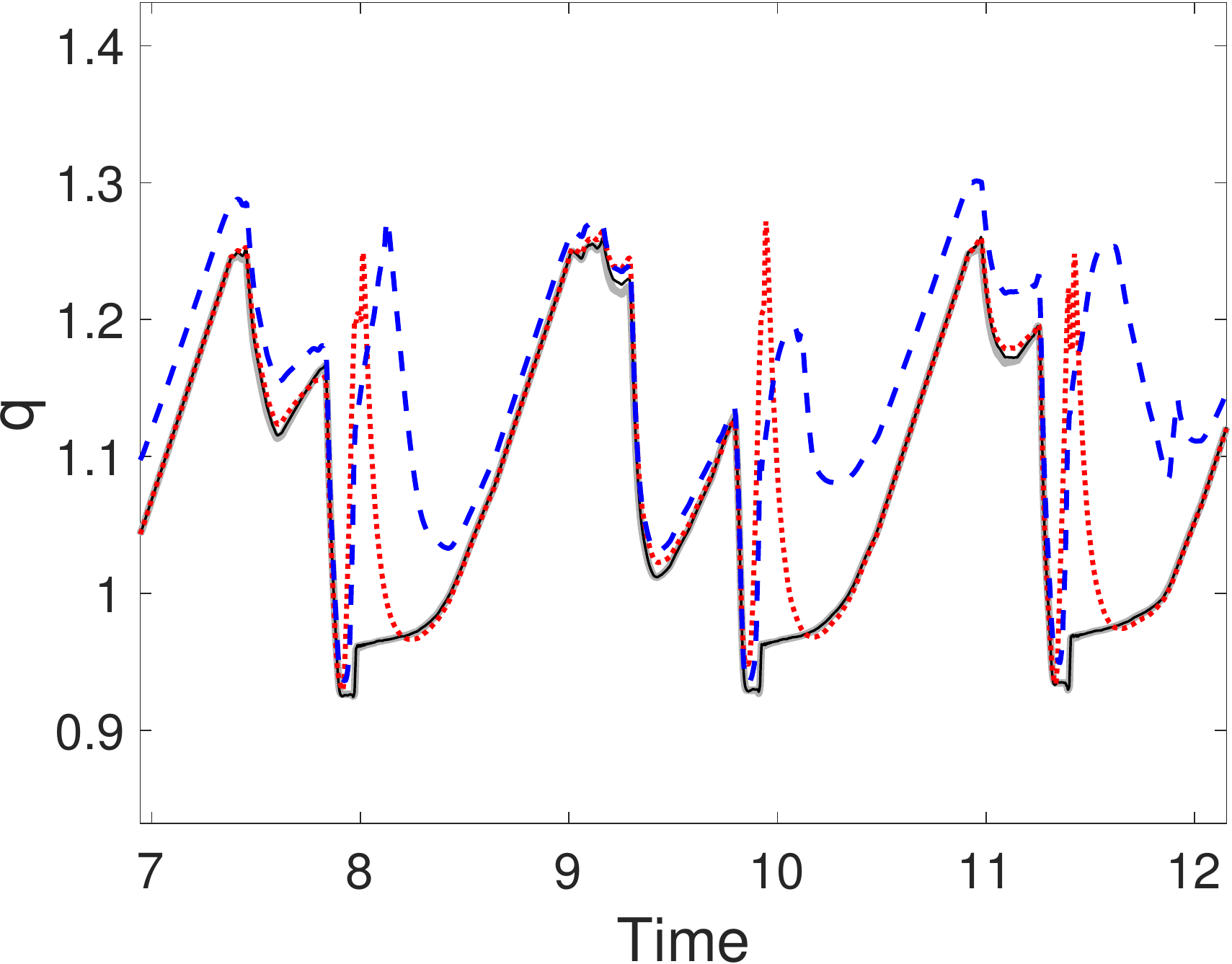}}
\subfigure[]{\includegraphics[width=0.24\columnwidth]{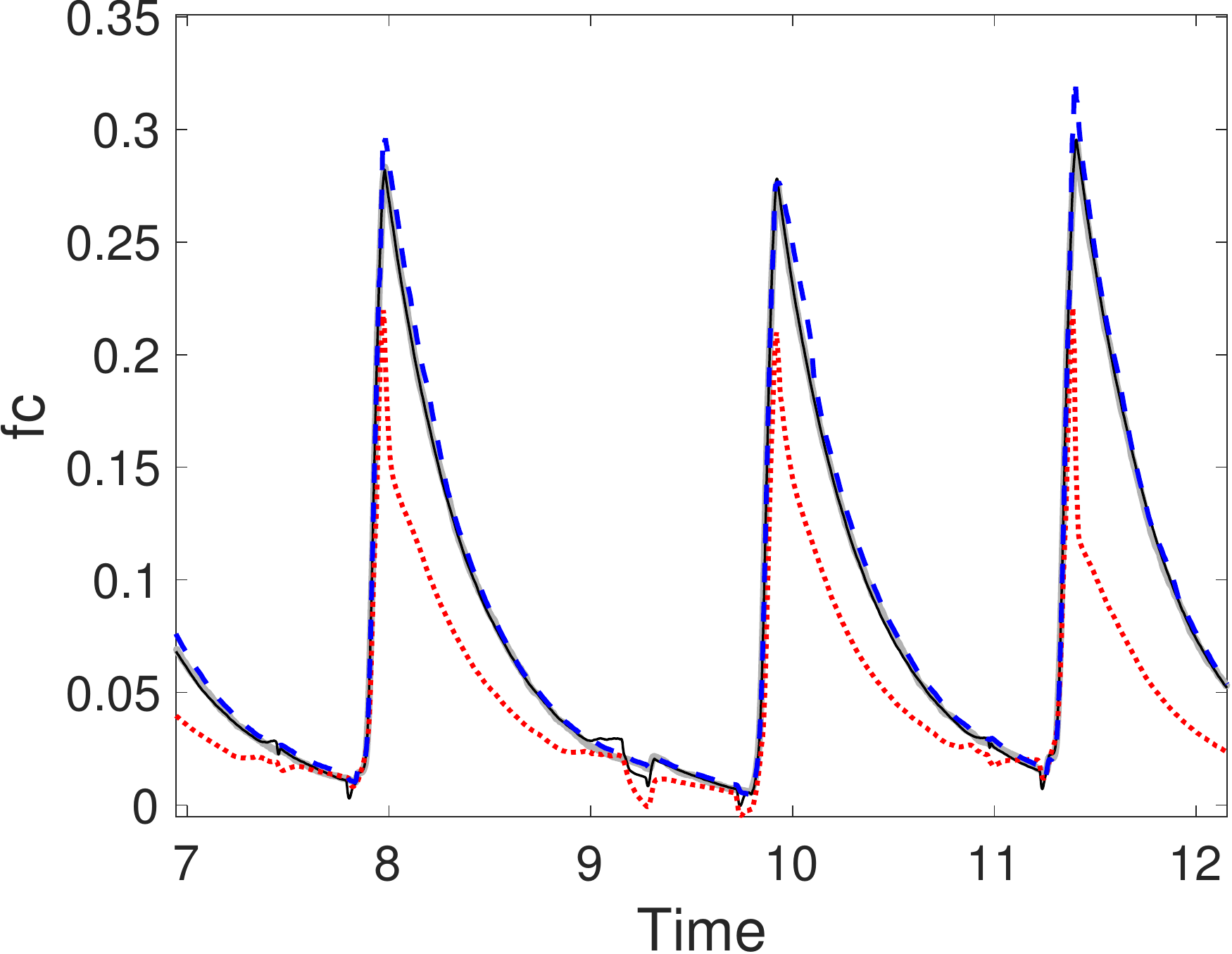}}
\subfigure[]{\includegraphics[width=0.24\columnwidth]{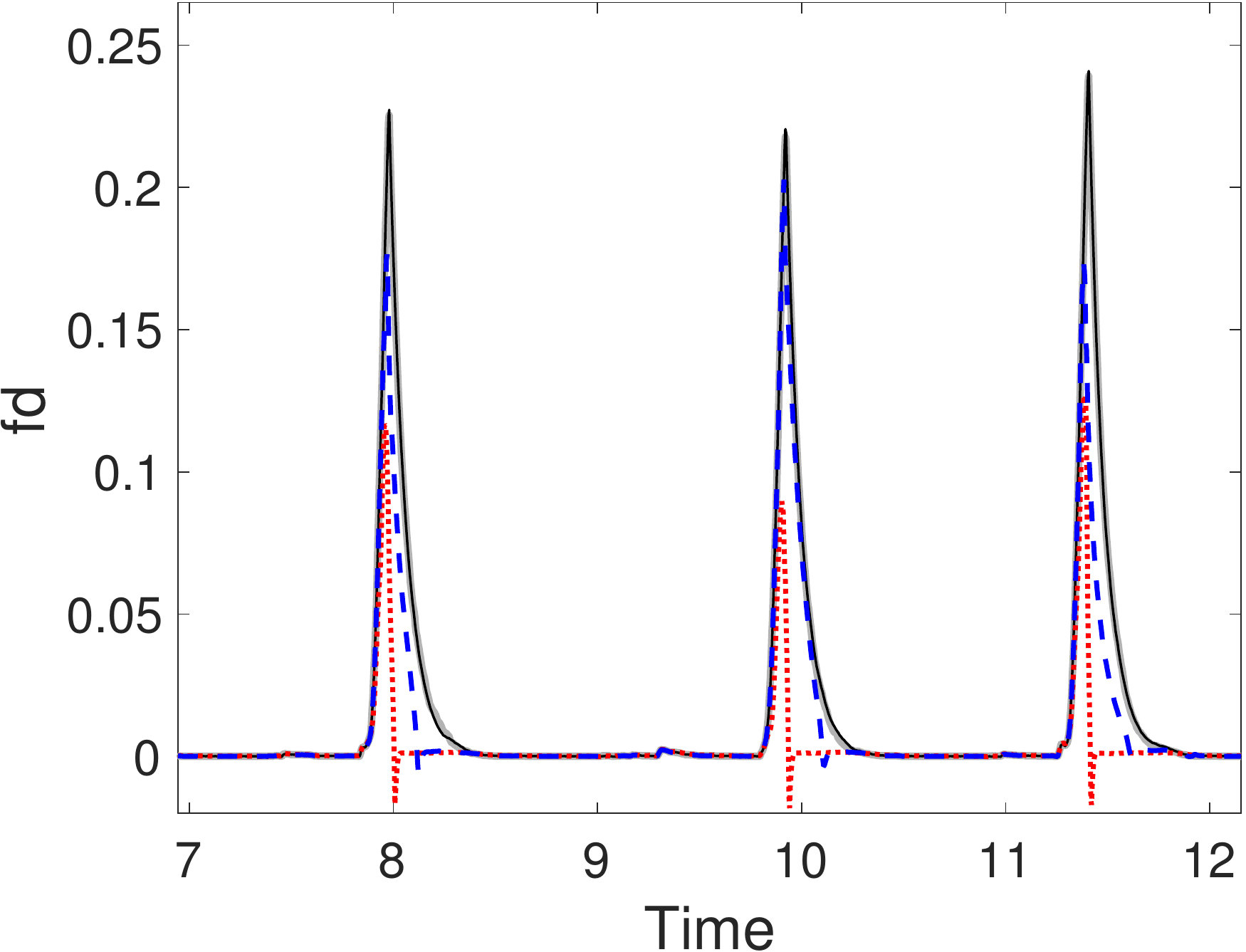}}
\subfigure[]{\includegraphics[width=0.24\columnwidth]{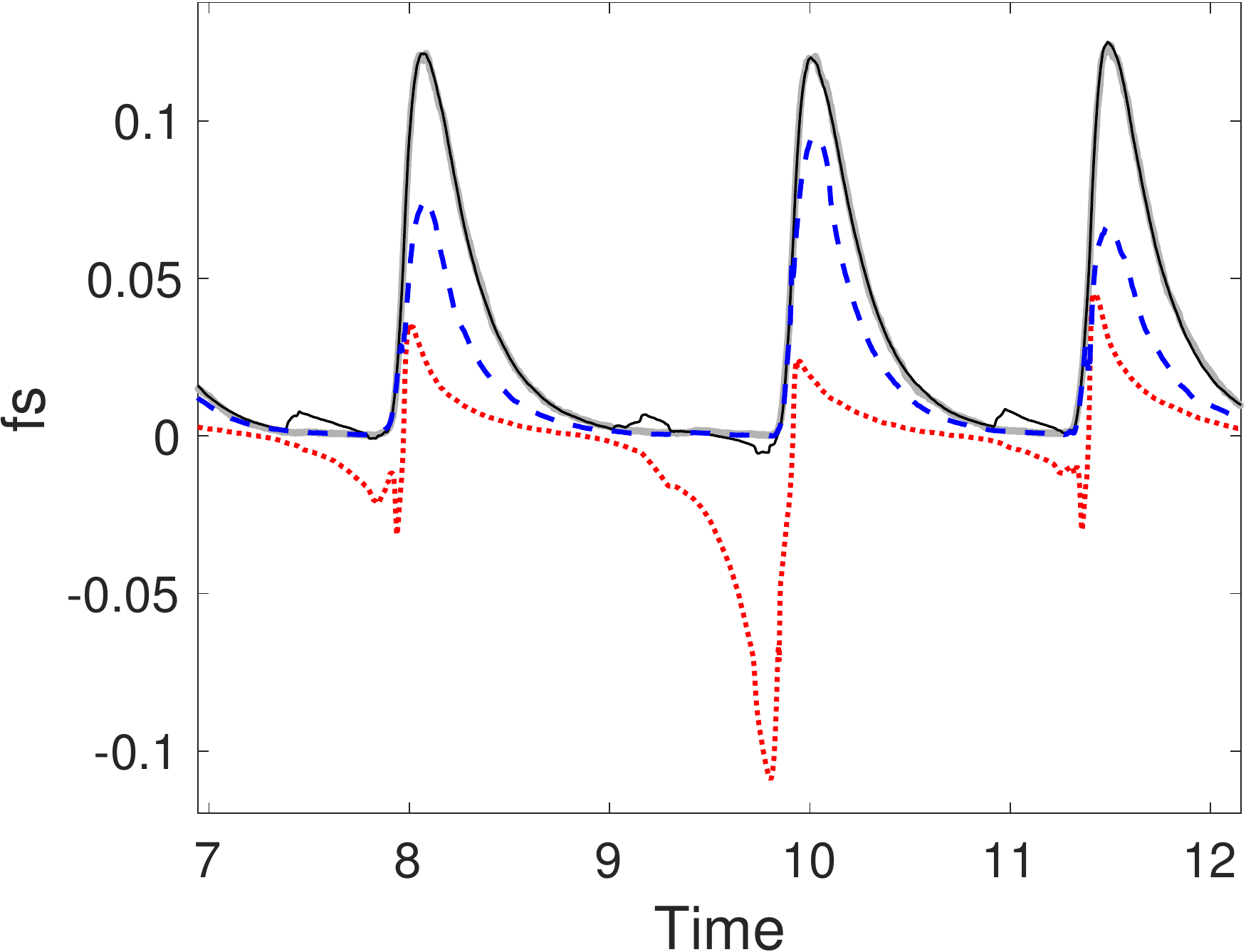}}
\caption{(a) True observations (red, dashed) incorporating cloud information are compared to the incorrect observation function (black, solid) which sets all the cloud fractions to zero in the RTM. (b-h) True state (gray, thick curve) compared to the result of filtering with the true observation function (black), the wrong observation function using only inflation of the observation covariance matrix (red, dashed) and the wrong observation function with iterative observation model error correction (blue, dashed).}
\label{RTM}
\end{figure}

\begin{figure}[h]
\center
\subfigure[]{\includegraphics[width=0.24\columnwidth]{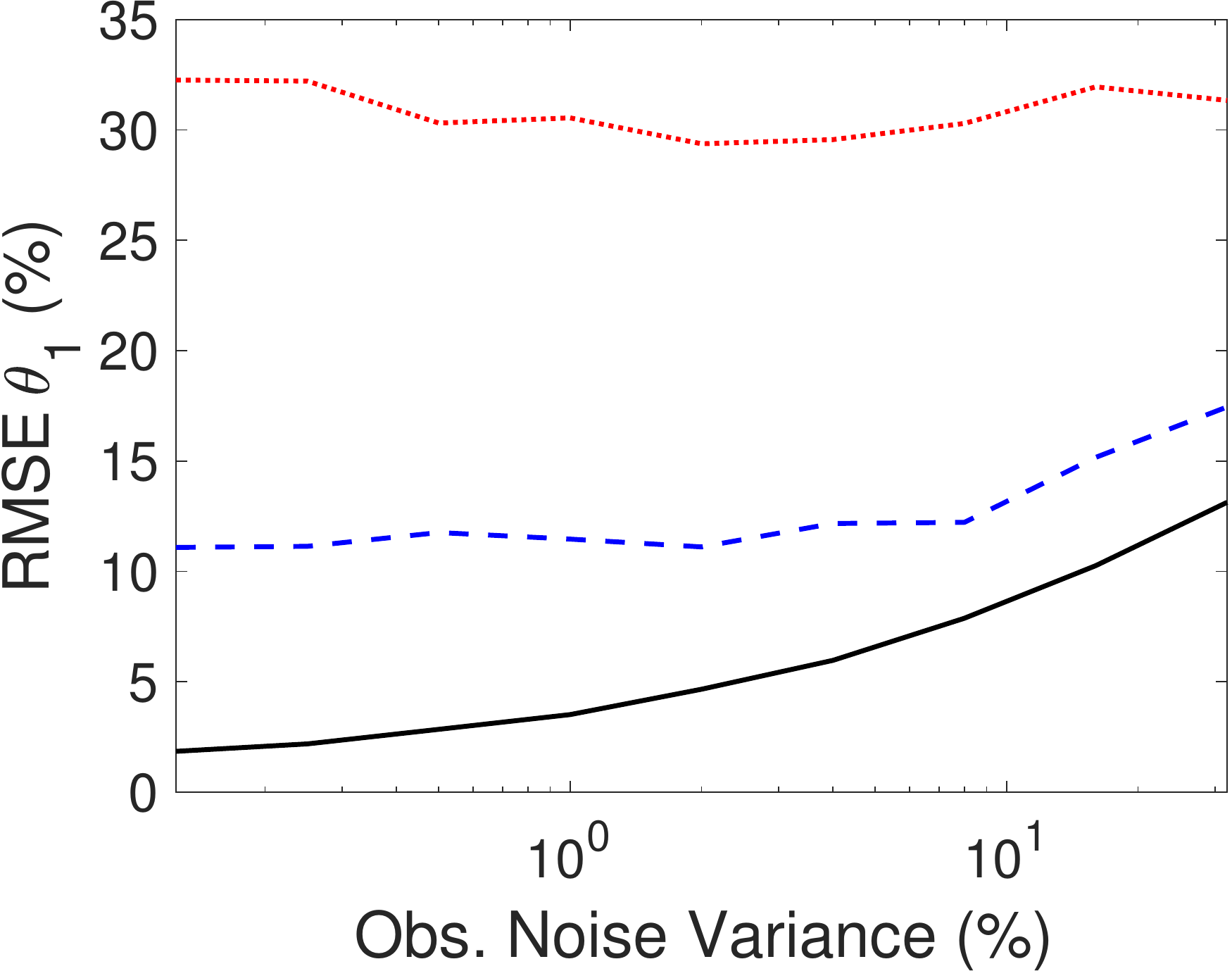}}
\subfigure[]{\includegraphics[width=0.24\columnwidth]{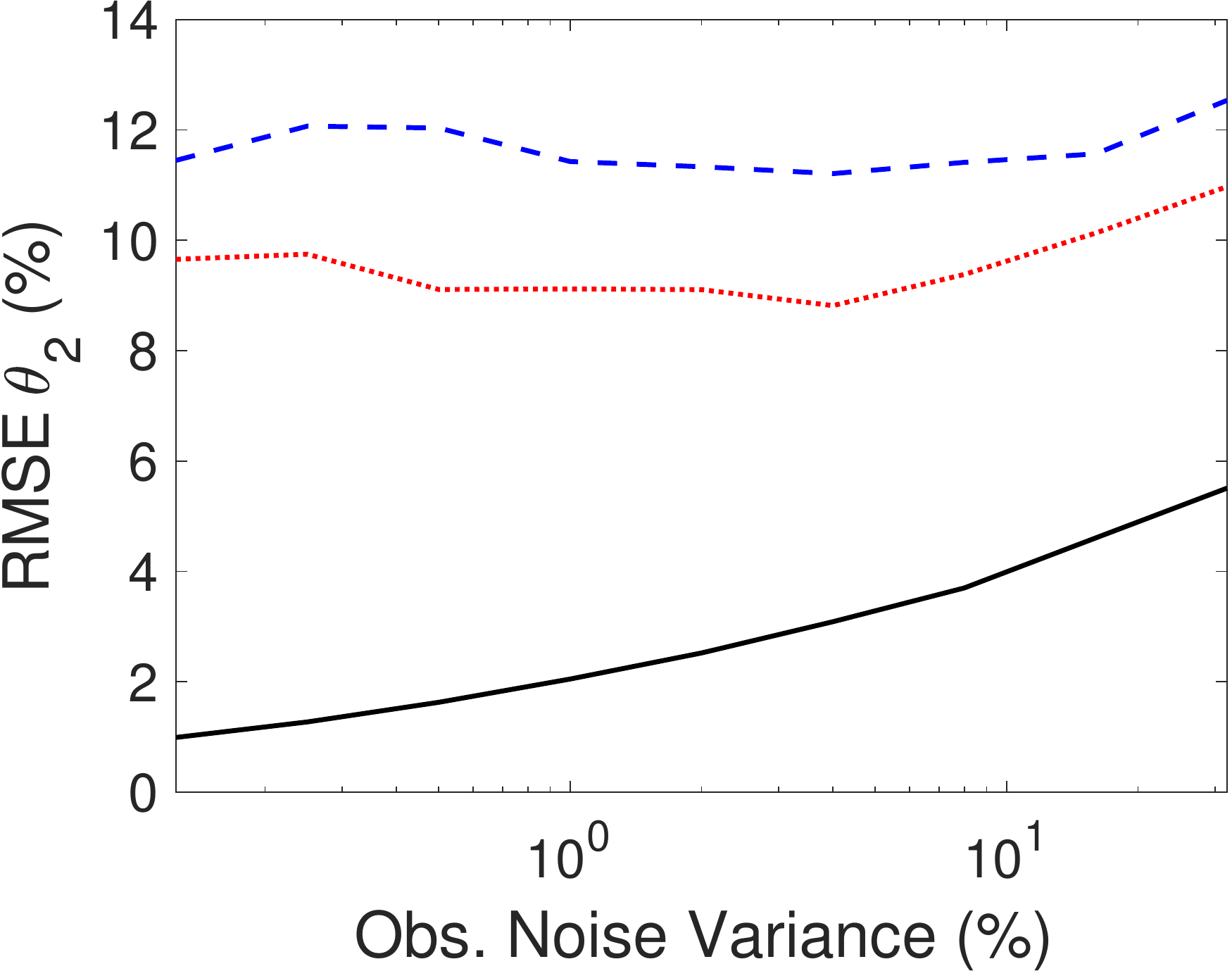}}
\subfigure[]{\includegraphics[width=0.24\columnwidth]{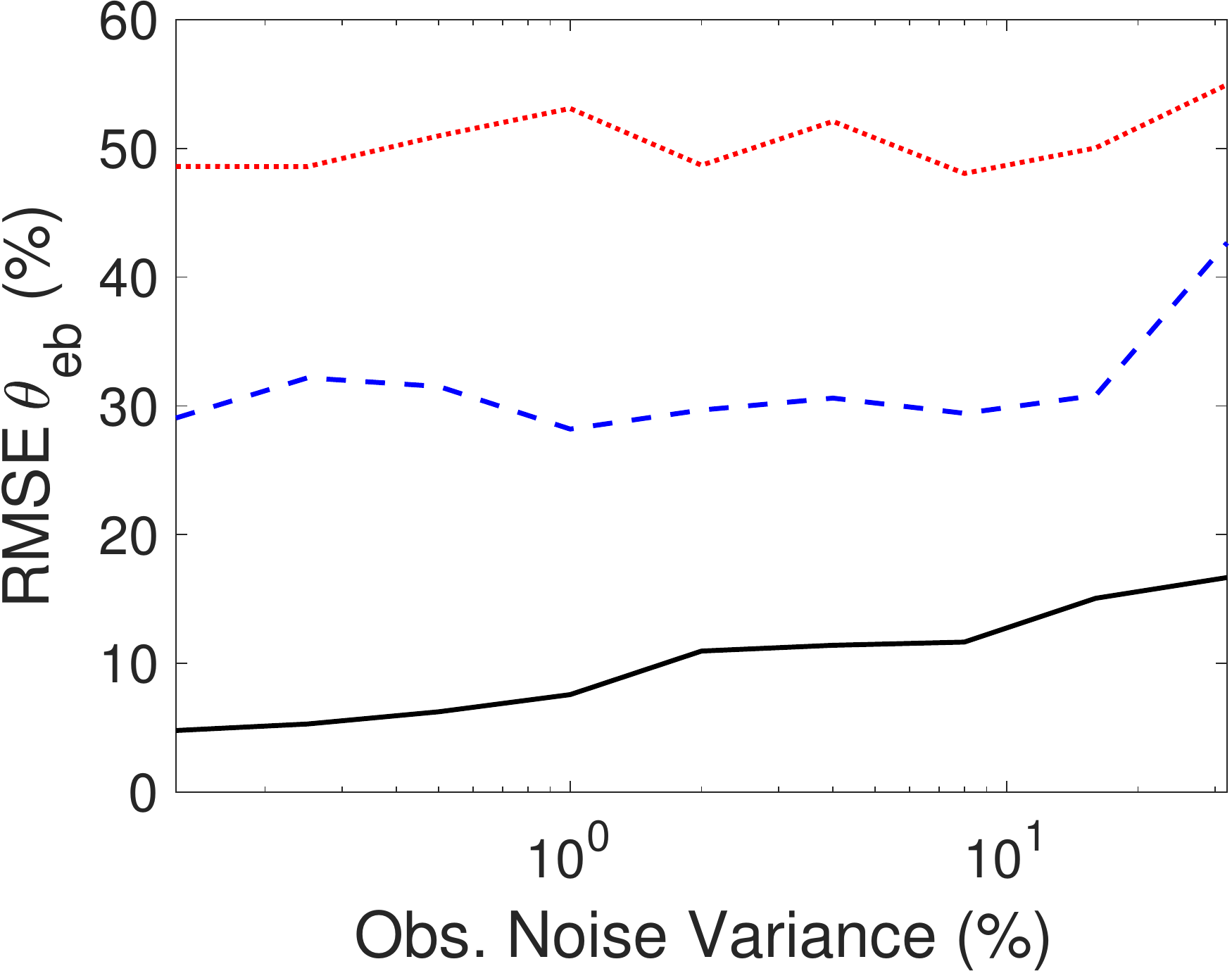}}
\subfigure[]{\includegraphics[width=0.24\columnwidth]{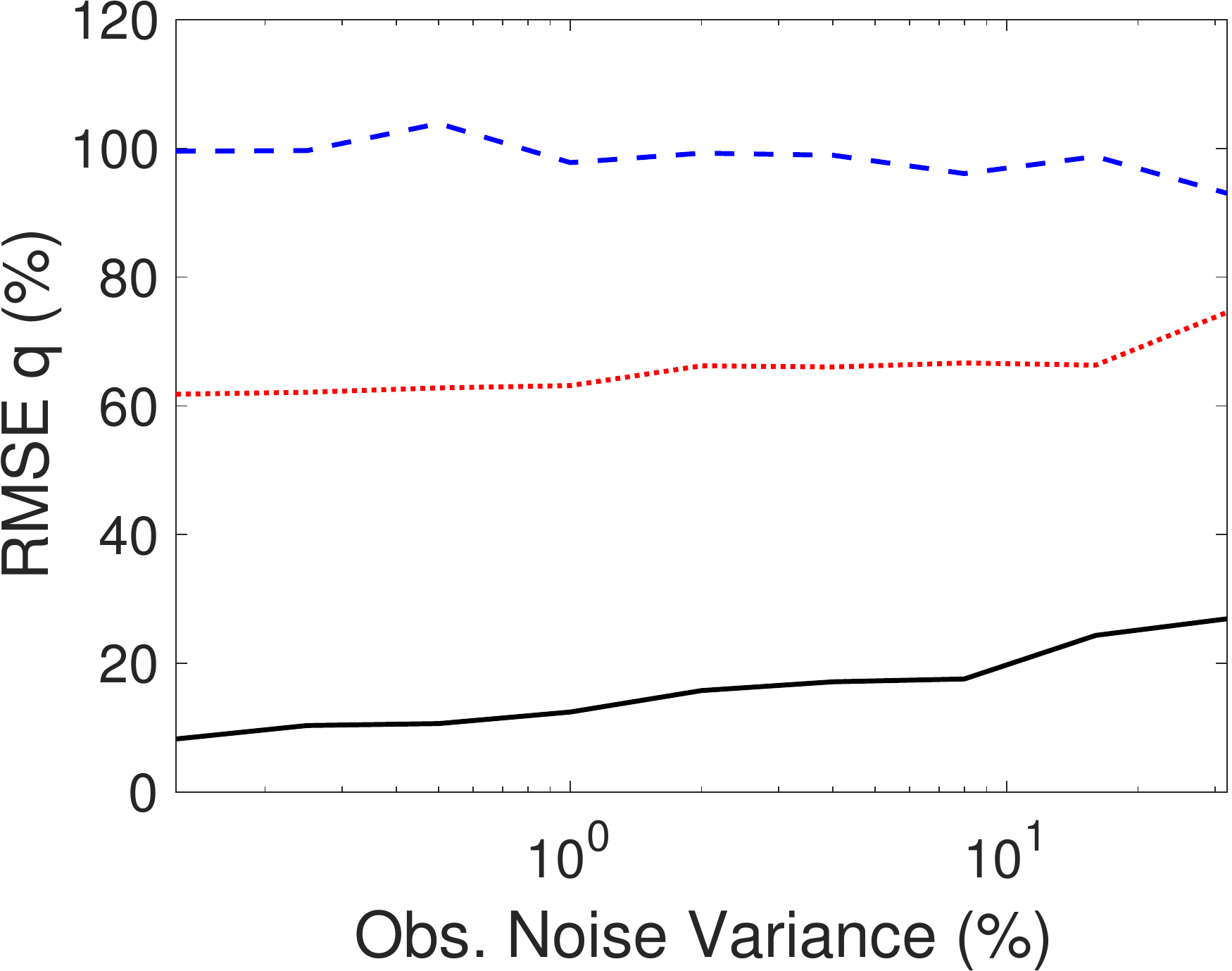}}
\subfigure[]{\includegraphics[width=0.24\columnwidth]{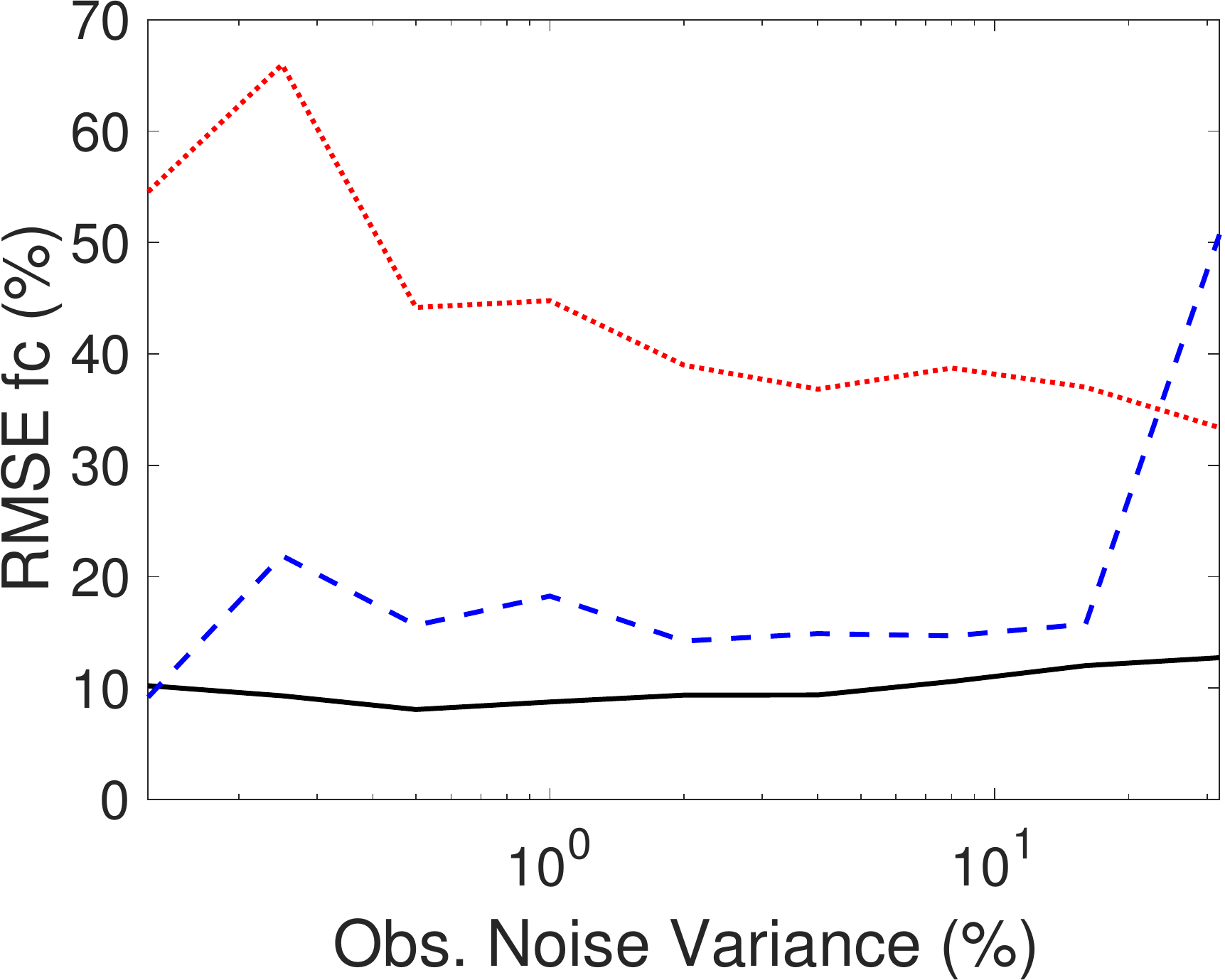}}
\subfigure[]{\includegraphics[width=0.24\columnwidth]{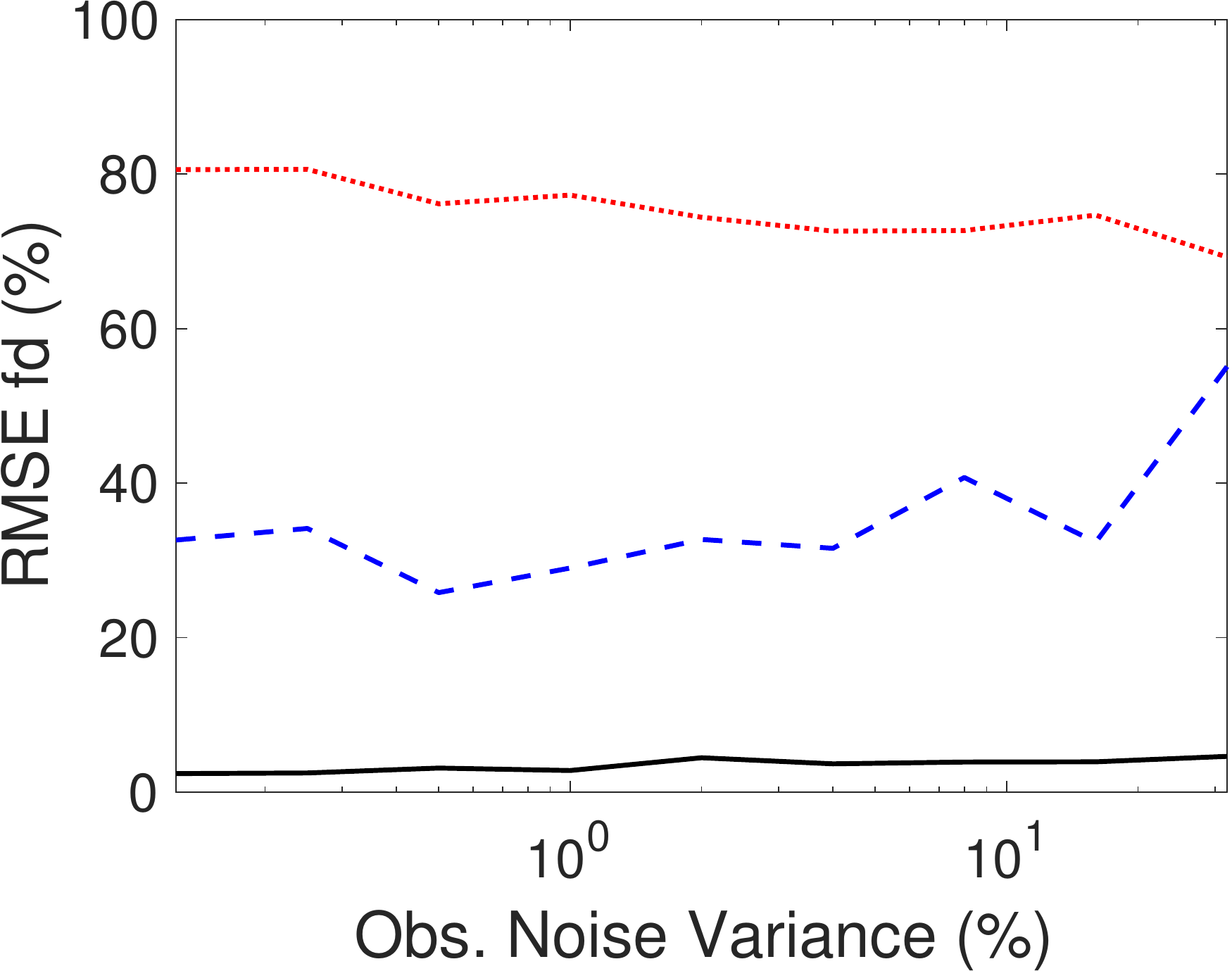}}
\subfigure[]{\includegraphics[width=0.24\columnwidth]{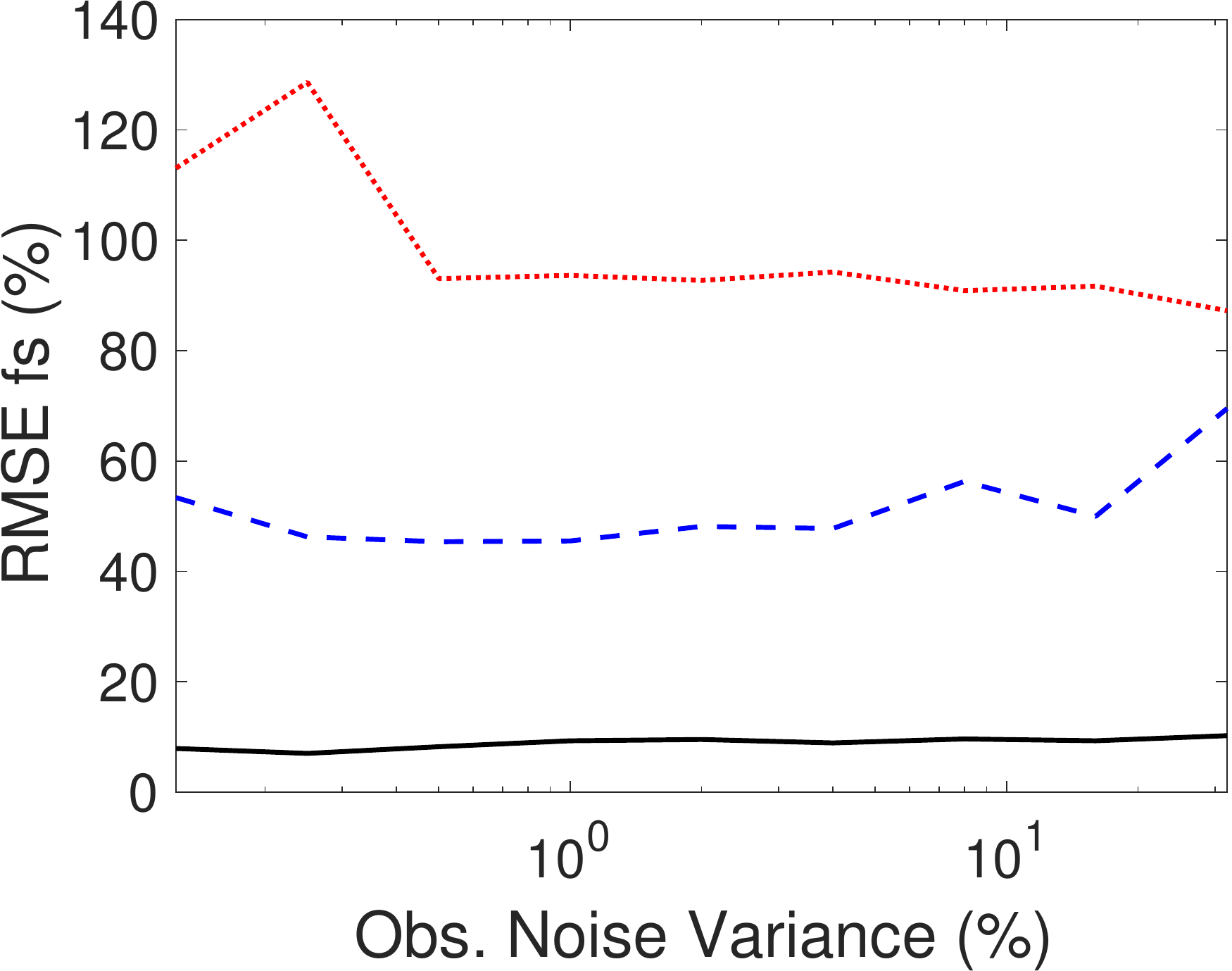}}
\caption{Robustness of filter estimates.  RMSE as a percentage of the standard deviation of each variable is shown as a function of observation noise percentage (noise variance is the given percentage of the the observation variance for each observed variable).  The filter using the true observation function (black, solid curve) is compared to the result of filtering with the wrong observation function using only inflation of the observation covariance matrix (red, dashed) and the wrong observation function with iterative observation model error correction (blue, dashed).}
\label{RTM2}
\end{figure}

The model  \cite{kbm:10} presented here represents a single column of atmosphere with three temperature variables $\theta_1,\theta_2,\theta_{eb}$ and a vertically averaged water vapor variable $q$.  The RTM also contains a stochastic multicloud parameterization with three variables $f_c,f_d,$ and $f_s$ which represent fractions of congestus, deep, and stratiform clouds respectively.  The three temperature variables are extrapolated to yield the temperature as a continuous function of the height, and then a simplified RTM can be used to integrate over this vertical profile to determine the radiation at various frequencies (see Berry and Harlim \cite{berry2017} for details).  We follow Liou \cite{liou:02} to incorporate information from the cloud fractions into the RTM in order to produce synthetic `true' observations at 16 different frequencies.  Each frequency has a different height profile which is integrated against the vertical temperature profile.  The presence of the different types of clouds influences these height profiles to simulate the cloud `blocking' radiation from below it.  We first show that the EnKF is capable of recovering most of the state variables from the observations when the correct observation model is specified (meaning the RTM includes the cloud fraction information from the model).  In Fig.~\ref{RTM} we show the true state (grey) along with the estimates produced using the correct observation model (black).

Next, we assume that the cloud fractions are unknown or that their effect on the RTM is poorly understood, and we attempt to assimilate the true observations using an RTM where the cloud fractions are held constant at zero (note that the cloud fractions are still present and evolving in the model used by the filter, but they are not included in the RTM used for the observation function of the filter).  We should note that this assimilation is impossible without artificially inflating the observation covariance matrix $R$ by a factor of 100.  The results of assimilating are shown in Fig.~\ref{RTM} (red, dotted).  Finally we apply the iterative observation model error correction (3 iterations) and the results are shown in Fig.~\ref{RTM} (blue, dashed).  Similar to the results of Berry and Harlim \cite{berry2017} the water vapor variable, $q$ is difficult to reconstruct in the presence of observation model error, however the cloud and temperature variables are significantly improved.  

\medskip

\begin{tabular}{| l | c | c | c | c | c | c | c |}
\hline
Percent Error (RMSE) & $\theta_1$ & $\theta_2$ & $\theta_{eb}$ & $q$ & $f_c$ & $f_d$ & $f_s$ \\
\hline
True Observation Function  & 2.8   & 1.6 &   6.2  & 10.6  &  8.1 &  3.1  &  8.2  \\
Wrong Observation Function &  30.3  &   9.1 & 51.0  & 62.8 &  44.2 &  76.2 &  93.1 \\
Model Error Correction & 11.8 &  12.0 &  31.5 & 103.9  & 15.6 &  25.8 &  45.4 \\
\hline
\end{tabular}

\medskip

In the table above we summarize the RMSE of each variable averaged over 4500 discrete filter steps (15.6 model time units with $dt=.0035$) for each filter, the observation noise variance was set at 0.5\% of the variance of each observed variable.  The observation model error correction is able to improve the estimation of all of the cloud fraction variables $f_c,f_d,$ and $f_s$ along with two of the temperature variables.  The estimation of $\theta_2$ was only slightly degraded.  The estimation of $q$ was more significantly degraded by the observation model error correction, probably because $q$ does not enter into the observation function as directly as the other variables.  These results compare favorably with Berry and Harlim \cite{berry2017} who also found that the $q$ variable was difficult to reconstruct in the presence of this observation model error, even using training data that included the true state.  

Since our approach here does not depend on perfect training data, we also found that our results were more robust to observation noise than the results of Berry and Harlim \cite{berry2017}.  In that approach, this was a significant issue  since it was assumed that the observation noise was small in order to be able to recover the true model error from the training data. As a result, the results were only robust up to observation noise levels of about 1\% of the variance of the observations. 

 In Fig.~\ref{RTM2} we show the robustness of the observation model error correction proposed here to increasing levels of observation noise. We find that the iterative observation model error correction is robust at noise levels over 10\% of the variance of the observations. At extremely low noise levels, such as levels near 0.1\%, the method of \cite{berry2017} has performance comparable to the true observation function, so when perfect full state training data is available and observation noise is small the methods have roughly equivalent behavior.

\section{Discussion}

Accurate linear and nonlinear filtering depends on thorough knowledge of model dynamics and the function connecting states to observations. The method proposed here uses an alternating minimization approach to iteratively correct observation model error, assuming knowledge of the correct dynamical model. This approach was shown to succeed in temporal and spatiotemporal examples as well as a cloud model. 

Although the iteration converges to eliminate observation model error in a wide variety of examples, there is no proof of global convergence of the method. This is typical for alternating minimization methods. A better understanding of the basin of convergence would be helpful, and the object of further study. 

The increasing diversity of measurement devices used in meteorological data assimilation is subject to a wide variety of separate errors. It is possible that more refined versions of the method can be designed to target particular subsets of the total observation error. The proof of concept carried out in this article show the potential for a relatively simple iterative solution to the problem, that can result in significant improvement in total RMSE.

We envision additional applications in other science and engineering areas, including hydrology, physical and biological experiments. A particular problem of interest in physiology is the common usage of intracellular neural models to assimilate extracellular measurements from single electrodes and electrode arrays. The observation function that connects such measurements to the model is not well understood by first principles and may vary by preparation. An automated way to solve this issue would potentially be a significant advance in data assimilation for neuroscience problems.

\bibliography{biasBib}

\section{Acknowledgements}

This research was partially supported by grants RTG/DMS-1246991 and DMS-1723175 from the National Science Foundation.

\end{document}